\documentclass[12pt]{amsart}
\usepackage{amsmath,amssymb,latexsym,cite,mathrsfs}
\usepackage{verbatim,wasysym,cite}
\usepackage[left=2.8cm,right=2.8cm,top=2.9cm,bottom=2.9cm]{geometry}
\usepackage{xcolor}

\usepackage[utf8]{inputenc}
\usepackage{microtype}
\usepackage{color,enumitem,graphicx}
\usepackage[colorlinks=true,urlcolor=blue, citecolor=red,linkcolor=blue,
linktocpage,pdfpagelabels, bookmarksnumbered,bookmarksopen]{hyperref}

\usepackage[english]{babel}

\renewcommand{\epsilon}{{\varepsilon}}

\numberwithin{equation}{section}
\newtheorem{theorem}{Theorem}[section]
\newtheorem{lemma}[theorem]{Lemma}
\newtheorem{remark}[theorem]{Remark}
\newtheorem{definition}[theorem]{Definition}
\newtheorem{proposition}[theorem]{Proposition}
\newtheorem{corollary}[theorem]{Corollary}

\DeclareMathOperator{\dom}{dom}
\DeclareMathOperator{\re}{Re}
\DeclareMathOperator{\ran}{ran}
\DeclareMathOperator{\im}{Im}
\DeclareMathOperator{\ess}{ess}
\DeclareMathOperator{\sech}{sech}
\DeclareMathOperator{\supp}{supp}
\DeclareMathOperator{\Span}{span}
\title[Nonlinear  Schr\"{o}dinger system on a star graph  ]
{Variational and stability properties of   coupled  NLS equations on the star graph}

\author[Liliana Cely and  Nataliia Golshchapova]{}
\email{mlcelyp@ufmg.br, nataliia@ime.usp.br }

\subjclass[2010]{Primary: 35Q55; Secondary: 35Q40}
\keywords{$\delta$ coupling, ground state, Nonlinear Schr\"{o}dinger equation,  orbital stability}

\begin{document}

\maketitle

\centerline{\scshape Liliana Cely and  Nataliia Goloshchapova}
{\footnotesize
\centerline{Institute of Mathematics and Statistics, University of S\~ao  Paulo,}
\centerline{Rua do Mat\~ao, 1010,  S\~ao  Paulo-SP, 05508-090, Brazil}
}

\begin{abstract}
We consider variational and stability properties of a    system of two coupled nonlinear Schr\"{o}dinger equations on the star graph $\Gamma$ with the $\delta$ coupling at the vertex of $\Gamma$.  The first part is devoted to the proof of an existence of the ground state as the minimizer of the constrained energy in the cubic case. This result extends the one obtained recently for the coupled NLS equations on the line. 

In the second part, we study stability properties of several families of standing waves in the case of  a general power nonlinearity. In particular, we study one-component standing waves $e^{i\omega t}(\Phi_1(x), 0)$ and $e^{i\omega t}(0, \Phi_2(x))$.  Moreover, we study two-component standing waves $e^{i\omega t}(\Phi(x), \Phi(x))$ for the case of power nonlinearity depending on a unique power parameter $p$. 

 To our knowledge, these are the first  results on variational and stability properties  of coupled NLS equations on graphs.
\end{abstract}

\medskip

\section{Introduction}
The  nonlinear Schr\"{o}dinger equation with the  focusing power nonlinearity  and the  $\delta$ coupling on the  star graph $\Gamma$
\begin{equation}\label{Eq0} 
i\partial_{t}u(t,x)+\Delta_{\gamma}u(t,x)+\left|u(t,x)\right|^{q-2}u(t,x)=0
\end{equation}
has been extensively studied during the last decade (see \cite{NQTG1, AQFF2, AngGol18,   Kai19, Noja1, NataMasa2020}). 
Here $\gamma\in\mathbb{R}\setminus\{0\},\,\,q>2, \, u: \mathbb{R} \times \Gamma \rightarrow \mathbb{C}^{N}$, and $\Delta_{\gamma}$ is the Laplace operator on $L^{2}(\Gamma)$ with the  $\delta$ coupling: for $v=(v_e)_{e=1}^N$
\begin{equation*}\label{Ku}
\begin{split}
(-\Delta_{\gamma}v)(x)=\left(-v''_{e}\right(x))_{e=1}^N,\qquad\dom(-\Delta_{\gamma})=\left\{v\in H^{2}(\Gamma):\,\sum^{N}_{e=1}v^{\prime}_{e}(0)=-\gamma v_{1}(0)\right\}.
\end{split}
\end{equation*}
The  well-posedness of \eqref{Eq0} was established  in \cite{AQFF,ArdCel21,  NataMasa2020}, whereas  the existence and the stability of standing waves  were studied in \cite{AQFF, NQTG, AngGol18, AngGol18a,  Kai19}.

Nonlinear PDEs
on graphs appear as mathematical models which describe various physical phenomena. In particular, \eqref{Eq0}  appears as a preferred model in optics of nonlinear Kerr media and dynamics of Bose-Einstein condensates (see \cite{BecImp15, QGAH, JolSem11, Noja1}). Graph models  arise as an approximation of  multidimensional narrow waveguides when their thickness parameters converge to zero (see \cite{Exner, Gnutz, Kuc02, Miros, Tokuno}). 

In the first part of this paper we  consider two  coupled nonlinear Schr\"{o}dinger equations
\begin{equation}\label{Eq}
\begin{cases} 
i\partial_{t}u(t,x)+\Delta_{\gamma}u(t,x)+\left(a\left|u(t,x)\right|^{2}+b\left|v(t,x)\right|^{2}\right)u(t,x)=0\\
i\partial_{t}v(t,x)+\Delta_{\gamma}v(t,x)+\left(b\left|u(t,x)\right|^{2}+c\left|v(t,x)\right|^{2}\right)v(t,x)=0\\
\left(u(0,x),v(0,x)\right)=\left(u_{0}(x),v_{0}(x)\right),
\end{cases} 
\end{equation}
where $\gamma>0$, $a,b,c\in\mathbb{R}$, $(t,x)\in \mathbb{R}\times\Gamma$,  $u,v:\mathbb{R}\times\Gamma\rightarrow\mathbb{C}^{N}$.
Coupled NLS equations
appear in many physical applications:   interaction of waves with different polarizations,  description of nonlinear modulations of two monochromatic waves,  interaction of Bloch-wave packets in a periodic system, 
evolution of two orthogonal pulse envelopes in birefringent optical fiber, evolution of two surface wave packets in deep water, the   Hartree-Fock theory for a double condensate (for the references see \cite{AkhAnk97, Bhat15, KivAgr03, Man74, Men87, WadIiz92}).

From the physical and the mathematical point of view, an interesting issue is to study the existence and the stability of standing waves of
the system \eqref{Eq}. For the principal results on the existence of solitary waves for the coupled NLS equations on $\mathbb{R}$ and $\mathbb{R}^N$ and their   stability/instability properties,  the reader is addressed to \cite{Barrst, Cipotti, Lop06, Lop11, Mai06,  Mai10, PelYan05,  Song,  NguWa, Wei}. A standing wave solution of \eqref{Eq} is a solution of the form $\left(e^{i\omega_{1}t}\Phi_1(x),e^{i\omega_{2}t}\Phi_2(x)\right)$, where $\omega_{1},\omega_{2}\in\mathbb{R}$ and $\left(\Phi_1,\Phi_2\right)$ solves  the following stationary problem
\begin{equation}\label{Eomega}
\begin{cases} 
-\Delta_{\gamma}\Phi_1 +\omega_{1}\Phi_1-\left(a\left|\Phi_1\right|^{2}+b\left|\Phi_2\right|^{2}\right)\Phi_1=0\\
-\Delta_{\gamma}\Phi_2 +\omega_{2}\Phi_2-\left(c\left|\Phi_2\right|^{2}+b\left|\Phi_1\right|^{2}\right)\Phi_2=0.
\end{cases} 
\end{equation}
It is classical idea  to look for the profile of the solitary wave of a Hamiltonian system as a solution of a certain minimization problem. In Section \ref{sec_2} we study the existence of the profiles $(\Phi_1(x), \Phi_2(x))$  being minimizers of the energy under the fixed mass constraint (depending on the constants $a,b,c, \gamma$). In particular, we find explicitly the minimizer $(e^{i\theta_1}\Phi_1(x), e^{i\theta_2}\Phi_2(x)),\, \theta_1, \theta_2\in\mathbb{R}$, where   $(\Phi_1, \Phi_2)$ is the solution to the stationary problem 
$$\begin{cases} -\Delta_\gamma\Phi_1+\omega \Phi_1-\frac{b^2-ac}{b-c}|\Phi_1|^2\Phi_1=0\\ -\Delta_\gamma\Phi_2+\omega \Phi_2-\frac{b^2-ac}{b-a}|\Phi_2|^2\Phi_2=0.\end{cases}$$

We managed to adapt the method of  \cite{NguWa} (elaborated for the coupled NLS equations on the line) for the case of the star graph.
This method requires that the  constants $a,b,c$ satisfy one  of the following two assumptions:

$(A1)$\, $0<b<\min\left\{a,c\right\}$ or 

$(A2)$\, $a,c>0$,  $b>\max\left\{a,c\right\}$.

The main idea is to use the generalization (developed in \cite{NQTG} and generalized in \cite{CacFin17} for the case of a general starlike graph) of the concentration-compactness principle \cite{Lions84a, Lions84b} for the case of the star graph  and the technique of  symmetric rearrangements on the star graph introduced in \cite{AQFF} (which is used to prove the absence of runaway case).
Moreover, we exploit the existence end the explicit form of the minimizer of the constrained energy  for the unique NLS equation with the $\delta$ coupling on $\Gamma$ obtained in \cite{NQTG}. 
The orbital stability of the standing wave associated with the minimizer follows standardly  (see Subsection \ref{subsec_2.3}). 

Notice that the results by \cite{NguWa} were recently extended in \cite{Bhat15} for the case of the generalized  power nonlinearity $$F_{p,q,r}(u,v)=(a|u|^{q-2}u+b|v|^p|u|^{p-2}u, c|v|^{r-2}v+b|u|^p|v|^{p-2}v).$$
It seems much more difficult to extend the technique from \cite{Bhat15}  for  the case of the star graph. In particular,  it is not clear how to prove that the corresponding variational problem is subadditive. The main difficulty is the presence of the term $\frac{N}{2}$ in   Pólya–Szegő inequality \eqref{deriv_rearrag}.

In the second part of the paper (see Section \ref{sec_3}) we deal with the stability properties of the standing waves for the system of coupled NLS equations on $\Gamma$.
In particular, in Subsection \ref{subsec_3.1} we study
one-component (one-hump) standing waves $(e^{i\omega t}\Phi_1(x),0)$ and $(0, e^{i\omega t}\Phi_2(x))$ and the nonlinearity $F_{p,q,r}(u,v)$.
The profiles $\Phi_1$ and $\Phi_2$ satisfy the following stationary equations
\begin{equation}\label{phi_12}
-\Delta_{\gamma} \Phi_1+\omega \Phi_1-a|\Phi_1|^{q-2} \Phi_1=0, \quad -\Delta_{\gamma} \Phi_2+\omega \Phi_2-c|\Phi_2|^{r-2} \Phi_2=0.\end{equation}
The description of each equation in  \eqref{phi_12} was  obtained in \cite[Theorem 4]{AQFF}. Namely, solutions constitute a family of $\left[\frac{N-1}{2}\right]+1$ vector-functions (one symmetric and the rest asymmetric). To prove spectral instability results for the family  we use an abstract   result from \cite{Grill88} (which permits to estimate the number of unstable eigenvalues $\lambda>0$).
For $p,q > 3$, using $C^2$ regularity of the  mapping
data-solution (associated with the corresponding Cauchy problem) and applying the abstract result from \cite{HenPer82}, we have shown orbital instability for the profiles  $\Phi_1, \Phi_2$.  This abstract result
states the nonlinear instability of a fixed point of a nonlinear mapping having the linearization $L$ of spectral radius
$r(L) > 1$. To apply the approach by \cite{Grill88} we need to estimate the Morse index of two self-adjoint in $L^2(\Gamma)$
operators associated with the second derivative of the action functional. These estimates were obtained in \cite[Theorem 3.1]{Kai19}.

In Subsection \ref{subsec_3.2} we study two-component  (multi-hump)
standing waves $(e^{i\omega t}\Phi(x), e^{i\omega t}\Phi(x))$ in the case of the one-parametric power nonlinearity 
$$F_{p}(u,v)=(|u|^{2p-2}u+b|v|^p|u|^{p-2}u, |v|^{2p-2}v+b|u|^p|v|^{p-2}v).$$ 
 It is easily seen that  $\Phi(x)$ solves
\begin{equation*}
-\Delta_\gamma \Phi+\omega \Phi-(b+1) \Phi^{2p-1}=0.
\end{equation*}
As in the previous case, the solutions  form a similar family of  $\left[\frac{N-1}{2}\right]+1$ vector-functions $(\Phi^\gamma_k, \Phi^\gamma_k), k=0, \ldots, \left[\frac{N-1}{2}\right]$.  Again we use the results by \cite{Grill88, HenPer82} mentioned above to prove spectral and orbital instability results for $\gamma>0, k\geq 1$  and $\gamma<0, k\geq 0$.
It is worth noting that in the multi-hump  case one of the main technical difficulties is that the ``real" part $\tilde L^R$ of a self-adjoint operator  associated with the second derivative of the action functional  is not  diagonal  (as it was in the previous case).   To overcome this difficulty we diagonalize  the system $\tilde L^R\vec{h}=\lambda\vec{h},\,\ \vec{h}=(h_1,h_2)$ making    linear transformation $h_+=h_1+h_2,\, h_-=h_1-h_2.$ 

Separately, we prove the orbital stability result for $\gamma>0, k=0$ (that is, we consider the   symmetric profile $\Phi_0^\gamma$ with decaying on $\mathbb{R}^+$ components). We generalize the approach by \cite{Lop11}. The key argument is to use the analytic perturbation theory to count the number of positive eigenvalues of the Hessian  associated with the action functional at $(\omega, \omega).$

In   Subsection \ref{subsec_3.3} we consider  $F_p(u,v)$ for $p=2, b=1$. In this situation, the system obeys 2D rotation invariance. Using \cite[Stability and Instability  Theorem]{GrilSha90}  we prove spectral instability results and orbital stability result for a standing wave generated by $2D$ rotation (stability is due to the centralizer subgroup). 

\subsection{Notation and some useful facts}
A star graph $\Gamma$ is constructed by the union of $N\geq2$ infinite half-lines connected at a single vertex $\nu$. Each edge $I_{e}$, $e=1,\ldots,N$, can be regarded as $\mathbb{R}^{+}$, and the vertex $\nu$ is placed at the origin. Given a function $u:\Gamma\rightarrow\mathbb{C}^{N}$, $u=(u_{e})^{N}_{e=1}$, where $u_{e}:\mathbb{R}^{+}\rightarrow\mathbb{C}$ denotes the restriction of $u$ to $I_{e}$. In particular, the nonlinear term in \eqref{Eq} is defined componentwise: $$\left\{\begin{array}{c}\left(a\left|u\right|^{2}+b\left|v\right|^{2}\right)u=\Big(\left(a\left|u_{e}\right|^{2}+b\left|v_{e}\right|^{2}\right)u_{e}\Big)_{e=1}^N,\\ \left(b\left|u\right|^{2}+c\left|v\right|^{2}\right)v=\Big(\left(b\left|u_{e}\right|^{2}+c\left|v_{e}\right|^{2}\right)v_{e}\Big)_{e=1}^N.\,\,\end{array}\right.$$
We denote by $u_{e}(0)$ and $u'_{e}(0)$ the limits of $u_{e}(x)$ and  $u'_{e}(x)$  as $x\rightarrow 0^{+}$. We say that a function $u$ is \textit{continuous} on $\Gamma$ if every restriction $u_{e}$ is continuous on $I_{e}$ and $u_{1}(0)=\ldots=u_{N}(0)$. The space of continuous functions is denoted by $C(\Gamma)$.

The natural Hilbert space associated with  the Laplace operator $-\Delta_\gamma$ is $L^2(\Gamma)$, which is defined as $L^{2}(\Gamma)=\bigoplus^{N}_{e=1}L^{2}(\mathbb{R}^{+})$.
The inner product in $L^2(\Gamma)$ is given by
\begin{align*}
 (u, v)_{2}=\operatorname{Re} \sum_{e=1}^{N}(u_{e}, v_{e})_{L^{2}(\mathbb{R}^{+})}, \quad u=\left(u_{e}\right)_{e=1}^{N}, v=\left(v_{e}\right)_{e=1}^{N}.
\end{align*}
 Analogously, for $1\leq p\leq\infty$, we define the space $L^{p}(\Gamma)$ as the set of functions on $\Gamma$ whose components belong to $L^{p}(\mathbb{R}^{+})$, and the norm is defined by
\begin{align*}
\left\|u\right\|^p_p=\sum^{N}_{e=1}\|u_{e}\|^{p}_{L^p(\mathbb{R}^+)},\,\,\,p\neq \infty,\,\,\,\,\,\,\,\,\,\,\,\,\|u\|_{\infty}=\max_{1\leq e\leq N}\|u\|_{\infty}.
\end{align*}
Depending on the context,   $\|\cdot\|_{p}$  will denote  the norm in $L^{p}$ either  on the graph $\Gamma$ or on the line. The Sobolev spaces $H^{1}(\Gamma)$ and $H^{2}(\Gamma)$ are  defined as
\begin{align*}
H^{k}(\Gamma)=\left\{u\in C(\Gamma):\,\,u_{e}\in H^{1}(\mathbb{R}^{+}),\,\,\,e=1,\ldots,N\right\},\, k=1,2.
\end{align*}
 The proof of the following proposition is a direct consequence of the Gagliardo-Nirenberg inequality on $\mathbb{R}^{+}$ (see \cite[I.31]{MiPec}).
\begin{proposition} Let  $q\in\left[2,\infty\right]$,  $1\leq p\leq q$, and $\mu=\frac{\frac{1}{p}-\frac{1}{q}}{\frac{1}{2}+\frac{1}{p}}$, then there exists $C>0$ such that
\begin{equation}\label{GNIq}
\left\|u\right\|_{q}\leq C\left\|u'\right\|^{\mu}_{2}\left\|u\right\|^{1-\mu}_{p}
\end{equation}
for any $u\in H^{1}(\Gamma)$.
\end{proposition}
\begin{remark}
Observe that for a compact graph (that is, an abstract graph without infinite edges),  a weaker version of the Gagliardo-Nirenberg inequality holds: \begin{equation}\label{GN_comp}
\left\|u\right\|_{q}\leq C\left\|u\right\|^{\mu}_{H^1}\left\|u\right\|^{1-\mu}_{p}.
\end{equation}
\end{remark}

We define the space $X$ to be the Cartesian product $H^{1}(\Gamma)\times H^{1}(\Gamma)$ equipped with the norm $\left\|(u,v)\right\|_{X}^2=\left\|u\right\|_{H^1}^2+\left\|v\right\|_{H^1}^2$, and  $X^{*}$ stands for its dual. The corresponding duality product is denoted by $\langle\cdot, \cdot\rangle_{X^*\times X}$.
\iffalse
For $\gamma\in\mathbb{R}$ the Laplace operator $-\Delta_{\gamma}$  has a precise interpretation as the self-adjoint operator on $L^{2}(\Gamma)$ associated with the quadratic form $F_{\gamma}$ defined on $H^{1}(\Gamma)$,
\begin{align*}
{F}_{\gamma}[u]&=\|u'\|^{2}_2-\gamma\left|u_{1}(0)\right|^{2}.
\end{align*}
It might be easily shown that  $\sigma_{\ess}(-\Delta_{\gamma})=\left[0,\infty\right)$ and  $\sigma_{p}(-\Delta_{\gamma})=\left\{\begin{array}{cc}
     \emptyset, \quad \gamma\leq0,  \\
     -\frac{\gamma^{2}}{N^{2}},\quad \gamma>0. 
\end{array}\right. $
\fi
Throughout this paper we use   $C, C_{\varepsilon}, C_{\varepsilon, \alpha, \beta}$ to denote various positive constants whose actual value is not important and which may vary from line to line. By $n(L)$  we denote the number of negative eigenvalues of a linear operator $L$ (counting multiplicities). 

\section{Variational analysis: the case of the cubic nonlinearity}\label{sec_2}
\subsection{Statement of the main result}
Define  the energy functional
\begin{equation*}\label{Euv1}
\begin{split}
E(u,v)=\|u'\|_2^2+\|v'\|_2^2-\gamma\left(|u_1(0)|^2+|v_1(0)|^2\right)-\frac{1}{2}\left(a\|u\|^4_4+c\|u\|^4_4\right)-b\|uv\|_2^2.
\end{split}
\end{equation*}
and the mass functionals
\begin{equation}\label{masses}
Q(u)=\|u\|_2^2,\qquad Q(v)=\|v\|_2^2
\end{equation}
associated to \eqref{Eq}.
  \begin{proposition} For any $\left(u_{0},v_{0}\right)\in X$ there exists  a unique maximal solution $(u,v)\in C\left(\mathbb{R},X\right)\cap C^{1}\left(\mathbb{R},X^{*}\right)$ of \eqref{Eq} satisfying $\left(u(0,x),v(0,x)\right)=(u_{0},v_{0})$. 
Moreover,
\begin{equation*}\label{QE}
E\left(u(t),v(t)\right)=E(u_{0}, v_0),\,\,\,\,\,\,\,Q(u(t))=Q(u_{0}),\,\,\,\,\,\,\,\, Q(v(t))=Q(v_{0})\,\,\,\,\,\,\,\mbox{for all}\,\,\,\,t\in\mathbb{R}.
\end{equation*}
\end{proposition}
\noindent The proof of the local well-posedness follows analogously to  \cite[Theorem 4.10.1]{Caz03} (see also Section 2 in \cite{AQFF}). For the proof of the global well-posedness see Remark \ref{general_system}.
\iffalse To show the global existence, assume that there exists $T<\infty$ such that $\lim\limits_{t\rightarrow T}\left\|(u(t),v(t))\right\|_{X}=\infty$. From the Cauchy-Schwartz inequality, Gagliardo-Nirenberg inequality \eqref{GNIq}, and the conservation
laws, we conclude the existence of  a constant $C$ that depends on $\left(u_{0},v_{0}\right)$ only, such that for all $(u,v)\in X$ (see \eqref{Gagli1}-\eqref{Gagli2}  for more detail)
\begin{equation}\label{contra}
E(u_{0},v_{0})=E\left(u(t),v(t)\right)\geq\left(\left\|\partial_{x}u(t)\right\|^{2}_{2}+\left\|\partial_{x}v(t)\right\|^{2}_{2}\right)\left(1-C\left(\left\|\partial_{x}u(t)\right\|^{-1}_{2}+\left\|\partial_{x}v(t)\right\|^{-1}_{2}\right)\right).
\end{equation}
Thus, letting $\left\|\partial_{x}u(t)\right\|^{2}_{2}+\left\|\partial_{x}v(t)\right\|^{2}_{2}\rightarrow\infty$ when $t\rightarrow T$ leads to a contradiction with \eqref{contra}. \fi

For given  real constants $a,b,c$ satisfying either $(A1)$ or $(A2)$ (see Introduction) and fixed frequency  $\omega>\frac{\gamma^{2}}{N^{2}}$ we  set
\begin{equation*}\label{lamda}
\alpha(\omega)=\frac{2(b-c)}{b^{2}-ac}(N\sqrt{\omega}-\gamma),\qquad\beta(\omega)=\frac{2(b-a)}{b^{2}-ac}(N\sqrt{\omega}-\gamma).
\end{equation*}
It is obvious that  $\alpha(\omega)$  and $\beta(\omega)$ run the interval $(0, \infty)$ when  $\omega$ runs the interval $(\frac{\gamma^2}{N^2}, \infty)$.

We are interested in  ground states of \eqref{Eomega}. By \textit{a ground state} we mean a minimizer $\left(\Phi_1,\Phi_2\right)\in X$  of the  energy $E$  in $H^1(\Gamma)$ constrained to the manifold of the states with fixed masses $Q(\Phi_1)=\alpha(\omega)$ and $Q(\Phi_2)=\beta(\omega)$. That is, we  study the constrained variational problem
\begin{equation}\label{JOTA}
J(\alpha,\beta)=\inf\left\{E(u,v)\,:\,u,v\in H^{1}(\Gamma),\,\, \left\|u\right\|^{2}_{2}=\alpha(\omega),\,\,\left\|v\right\|^{2}_{2}=\beta(\omega) \right\}.
\end{equation}
In what follows we will use notation $\alpha, \beta$ instead of $\alpha(\omega), \beta(\omega)$. 
A \textit{minimizing sequence} for \eqref{JOTA} is a sequence $\left\{\left(u_{n},v_{n}\right)\right\}$ in $X$ such that $Q(u_{n})=\alpha$, $Q(v_{n})=\beta$ for all $n\in\mathbb{N}$, and $\lim\limits_{n\rightarrow\infty}E(u_{n},v_{n})=J(\alpha,\beta)$. We denote the set of nontrivial
minimizers  by
\begin{equation*}
\mathcal{G}(\alpha,\beta)=\left\{(u,v)\in X\,:\,J(\alpha,\beta)=E(u,v),\, \left\|u\right\|^{2}_{2}=\alpha,\,\,\left\|v\right\|^{2}_{2}=\beta \right\}.
\end{equation*}

Below we state the main result of the first part.
\begin{theorem}\label{Theo01}
Suppose that the numbers   $a,b,c$ satisfy either (A1) or (A2), $\omega>\frac{\gamma^{2}}{N^{2}}$, $\gamma>0$, and $N\sqrt{\omega}-\gamma\leq \frac{2\gamma}{N}$. Then the following statements hold.\\
$(i)$  Any minimizing sequence $\left\{\left(u_{n},v_{n}\right)\right\}$ for $J(\alpha,\beta)$ is relatively compact in $X$. That is, there exist a subsequence $\left\{\left(u_{n_{k}},v_{n_{k}}\right)\right\}$ and  $\left(\Phi_1,\Phi_2\right)\in X$ such that $\{\left(u_{n_{k}},v_{n_{k}}\right)\}$ converges to $\left(\Phi_1,\Phi_2\right)$ in $X$. Therefore, there exists a minimizer for problem \eqref{JOTA}, and hence 
$\mathcal{G}(\alpha,\beta)$ is non-empty.\\
$(ii)$ For each minimizing sequence $\left\{\left(u_{n},v_{n}\right)\right\}$ we have:
\begin{equation}\label{EeQq}
\lim\limits_{n\rightarrow\infty}\inf_{(\Phi_1,\Phi_2)\in \mathcal{G}(\alpha,\beta)}\left\|(u_{n},v_{n})-(\Phi_1,\Phi_2)\right\|_{X}=0.
\end{equation}
\end{theorem}
\begin{remark}\label{rem_minimizer}
The restriction  $N\sqrt{\omega}-\gamma\leq \frac{2\gamma}{N}$ has to be commented.

$(i)$ The method that we use requires an existence of a minimizer for  problems  \eqref{J12}. Due to \cite[Theorem 1]{NQTG} the problem 
\begin{equation}\label{min_AN}\inf\left\{\frac{1}{2}\left(\|u'\|_2^2-\gamma|u_{1}(0)|^2-\frac{1}{2}\|u\|_4^4\right): u\in H^1(\Gamma),\, \|u\|_2^2=m\right\} \end{equation} has a solution for a small mass, i.e.  $m\leq m^*=\frac{4\gamma}{N}$ (see formula (1.4) in \cite{NQTG}). In problems \eqref{J12} the term $\frac{1}{2}\|u\|_4^4$ from \eqref{min_AN} has to be substituted by $\frac{b^2-ac}{2(b-c)}\|u\|_4^4$ and $\frac{b^2-ac}{2(b-a)}\|u\|_4^4$.  Then in this case we will have that $m^*$ coincides with  $\frac{4\gamma}{N}\frac{(b-c)}{b^2-ac}$ and $\frac{4\gamma}{N}\frac{(b-a)}{b^2-ac}$ respectively. Finally, recalling equalities   $\alpha=\frac{2(b-c)}{b^2-ac}(N\sqrt{\omega}-\gamma)$, $\beta=\frac{2(b-a)}{b^2-ac}(N\sqrt{\omega}-\gamma)$, and assuming that $m$ coincides with $\alpha$ or $\beta$, we obtain that the restriction $m\leq m^*$ is equivalent to $N\sqrt{\omega}-\gamma\leq \frac{2\gamma}{N}$.

$(ii)$ It is interesting to compare the interval $(\frac{\gamma^2}{N^2},  \frac{1}{N^2}\left(\frac{2\gamma}{N}+\gamma\right)^2] $ of existence of solution to \eqref{JOTA}  with  the  one given by  \cite[Theorem 3]{CacFin17}. In particular, the theorem  states that for  $\omega\in (E_0, E_0+\delta)$  (with sufficiently small $\delta$) the stationary equation  
$$ H\Phi+\omega\Phi-|\Phi|^{q-2}\Phi=0$$
has a unique solution. Here $H$ is the Laplace operator with the $\delta$ coupling and  a linear potential on a general starlike graph, and $-E_0=\inf\sigma(H)$. 

In \cite[Remark 4.1]{AQFF} the authors mention that the restriction $m\leq m^*$ is not optimal, and the interval of existence of solution to  minimization problem \eqref{min_AN} is bigger.  This fact suggests that the solution to \eqref{JOTA} exists in some interval $(\frac{\gamma^2}{N^2}, \frac{\gamma^2}{N^2}+\delta)\supset (\frac{\gamma^2}{N^2}, \frac{1}{N^2}\left(\frac{2\gamma}{N}+\gamma\right)^2]$. 
\end{remark}
We have managed to obtain  the explicit characterization of the set of minimizers $\mathcal{G}(\alpha,\beta)$.
\begin{theorem}\label{theo02}
Suppose that assumptions of Theorem \ref{Theo01} hold. Then  the set of ground states
is given by
\begin{equation*}
\mathcal{G}(\alpha,\beta)=\left\{\left(e^{i\theta_1}\sqrt{\frac{b-c}{b^{2}-ac}}\Phi_{\omega,\gamma}(\cdot),e^{i\theta_2}\sqrt{\frac{b-a}{b^{2}-ac}}\Phi_{\omega,\gamma}(\cdot)\right)\,\,:\theta_1,\theta_2\in\mathbb{R}\right\},
\end{equation*}
where $(\Phi_{\omega,\gamma})_{e}=\phi_{\omega, \gamma}$, $e=1,\ldots,N$, with
\begin{equation}\label{Phi}
 \phi_{\omega, \gamma}(x)=\sqrt{2\omega}\,\sech\left(\sqrt{\omega}x+\tanh^{-1}\left(\frac{\gamma}{N\sqrt{\omega}}\right)\right).
\end{equation}
\end{theorem}
Next we give the definition of an orbital stability. 
\begin{definition}
The standing wave $\left(u(t,x),v(t,x)\right)=\left(e^{i\omega t}\Phi_1(x), e^{i\omega t}\Phi_2(x)\right)$ is said to be orbitally stable in $X$ if for any $\epsilon>0$ there exists $\eta>0$ with the following property: if $(u_{0},v_{0})\in X$ satisfies $$\left\|\left(u_{0},v_{0}\right)-\left(\Phi_1,\Phi_2\right)\right\|_{X}<\eta,$$ then the solution $(u(t,x),v(t,x))$ of \eqref{Eq} with $(u(0,\cdot),v(0,\cdot))=(u_{0},v_{0})$ satisfies
\begin{equation*}
\sup_{t\geq0}\inf_{\theta_1,\theta_2\in\mathbb{R}}\left\|\left(u(t,\cdot),v(t,\cdot)\right)-\left(e^{i\theta_1}\Phi_1,e^{i\theta_2}\Phi_2\right)\right\|_{X}<\epsilon.
\end{equation*}
Otherwise, the standing wave $\left(e^{i\omega t}\Phi_1(x), e^{i\omega t}\Phi_2(x)\right)$ is said to be orbitally unstable in $X$.
\end{definition}
Using the arguments from \cite{NQTG, AQFF, NguWa}, the compactness of the minimizing sequences (see Lemma \ref{lemma10} below), and the  uniqueness of the ground state (up to phase) proved  in Theorem \ref{theo02}, one arrives at the following result.
\begin{corollary}\label{kokoro1}
Suppose that assumptions of Theorem \ref{Theo01} hold. Then the standing wave solution 
\begin{equation*}
\left(e^{i\omega t}\sqrt{\frac{b-c}{b^{2}-ac}}\Phi_{\omega,\gamma},e^{i\omega t}\sqrt{\frac{b-a}{b^{2}-ac}}\Phi_{\omega,\gamma}\right),
\end{equation*}
where $\Phi_{\omega,\gamma}$ is defined by \eqref{Phi}, is orbitally stable in $X$.
\end{corollary}

\subsection{Existence of ground states}
In this subsection we  prove  Theorem \ref{Theo01}. We have divided the proof into several lemmas and propositions. Throughout this subsection,  we assume that the hypotheses of Theorem \ref{Theo01} hold.
\begin{proposition}
For any $\alpha,\beta>0$  we have  $-\infty<J(\alpha,\beta)<0$.
\end{proposition}
\begin{proof}
Firstly, we show $J(\alpha,\beta)<0$. Let $u\in H^{1}(\Gamma)$ be  such that $\left\|u\right\|^{2}_{2}=\alpha$ and we take $v=\sqrt{\frac{\beta}{\alpha}}u$, then $\left\|v\right\|^{2}_{2}=\beta$. Define $u_{r}(x)=\sqrt{r}u(rx)$ and $v_{r}(x)=\sqrt{r}v(rx)$ for  $r>0$. Since $\left\|u_{r}\right\|^{2}_{2}=\alpha$, $\left\|v_{r}\right\|^{2}_{2}=\beta$, and $\left\|u_{r}'\right\|^{2}_{2}=r^{2}\left\|u'\right\|^{2}_{2}$, then  we obtain

\begin{equation*}
\begin{split}
E(u_{r},v_{r})=r^{2}\left(\left\|u'\right\|^{2}_{2}+\left\|v'\right\|^{2}_{2}\right)-r\gamma\left(\left|u_{1}(0)\right|^{2}+\left|v_{1}(0)\right|^{2}\right)
-\frac{r}{2}\left(a\|u\|^{4}_4+c\|v\|^{4}_4\right)-rb\|uv\|^{2}_2.
\end{split}
\end{equation*}
Since $v=\sqrt{\frac{\beta}{\alpha}}u$, then we have
\begin{equation*}
\begin{split}
E(u_{r},v_{r})&=\left(1+\frac{\beta}{\alpha}\right)\left(r^{2}\left\|u'\right\|^{2}_{2}-r\gamma\left|u_{1}(0)\right|^{2}\right)-r\frac{(b^{2}-ac)(2b-a-c)}{2(b-c)^{2}}\left\|u\right\|^{4}_{4}\\
&\leq\left(1+\frac{\beta}{\alpha}\right)r^{2}\left\|u'\right\|^{2}_{2}-r\frac{(b^{2}-ac)(2b-a-c)}{2(b-c)^{2}}\left\|u\right\|^{4}_{4}.
\end{split}
\end{equation*}
By  $(b^{2}-ac)(2b-a-c)>0$, we can choose $r$ small enough to ensure that $E(u_{r},v_{r})<0$, hence $J(\alpha,\beta)<0$. Secondly, we prove that $J(\alpha,\beta)>-\infty$. Using  Gagliardo-Nirenberg inequality \eqref{GNIq}, the  Cauchy-Schwartz inequality, and the Young inequality, we have
\begin{equation}\label{Gagli1}
\left|u_{1}(0)\right|^{2}\leq\left\|u\right\|^{2}_{\infty}\leq C\left\|u'\right\|_{2}\left\|u\right\|_{2}\leq \varepsilon \|u'\|_2^2+C_\varepsilon \|u\|_2^2,
\end{equation}
\begin{equation}\label{Gagli3}
\int_{\Gamma}\left|u\right|^{2}\left|v\right|^{2}dx\leq\left\|u\right\|^{2}_{4}\left\|v\right\|^{2}_{4}\leq\frac{1}{2}\left(\left\|u\right\|^{4}_{4}+\left\|v\right\|^{4}_{4}\right),
\end{equation}
and
\begin{equation}\label{Gagli2}
\left\|u\right\|^{4}_{4}=C\left\|u'\right\|_{2}\left\|u\right\|^{3}_{2}\leq \epsilon \left\|u'\right\|^{2}_{2}+C_{\epsilon}\left\|u\right\|^{6}_{2}
\end{equation}
for all $u\in H^{1}(\Gamma)$ and any $\epsilon>0$. Therefore, by \eqref{Gagli1}, \eqref{Gagli3},  and  \eqref{Gagli2}, we have
\begin{equation*}
\begin{split}
&E(u,v)=\|u'\|_2^2+\|v'\|_2^2-\gamma\left(|u_{1}(0)|^2+|v_{1}(0)|^2\right)-\frac{1}{2}\left(a\left\|u\right\|^{4}_{4}+c\left\|v\right\|^{4}_{4}\right)-b\|uv\|^{2}_2\\
&\geq \|u'\|_2^2+\|v'\|_2^2-\gamma\left(|u_{1}(0)|^2+|v_{1}(0)|^2\right)-C\left(\left\|u\right\|^{4}_{4}+\left\|v\right\|^{4}_{4}\right)\\&\geq (1-\varepsilon)\left(\|u'\|_2^2+\|v'\|_2^2\right)-C_\varepsilon\left(\|u\|_2^2+\|v\|_2^2+ \|u\|_2^6+\|v\|_2^6\right)\\&= (1-\varepsilon)\left(\|u'\|^2_{H^{1}}+\|v'\|^2_{H^{1}}\right)-(1-\varepsilon)(\alpha+\beta)-C_\varepsilon\left(\alpha+\beta+ \alpha^3+\beta^3\right)\\
&=(1-\epsilon)\left(\left\|u\right\|^{2}_{H^{1}}+\left\|v\right\|^{2}_{H^{1}}\right)-C_{\epsilon,\alpha,\beta}>-\infty.
\end{split}
\end{equation*}
This ends the proof.
\end{proof}
Next we introduce the concentration-compactness technique for $\Gamma$.
Let $x\in I_{e}$ and $y\in I_{j}$, $e,j=1,\ldots,N$, be  two points of graph. We define the distance   
\begin{equation}\label{55}
d(x,y):=
\begin{cases} 
\left|x-y\right|\,\,\,\mbox{for}\,\,\,\,e= j,\\
x+y\,\,\,\mbox{for}\,\,\,\,e\neq j,
\end{cases}
\end{equation}
and  the open ball of center $x$ and radius $r$
\begin{equation*}
B(x,r):=\left\{y\in \Gamma\,:\,\,d(x,y)<r\right\}.
\end{equation*}
 Let   $x\in I_{e}$,  we set  the $L^{p}$ norm restricted to the ball $B(x,r)$ 
\begin{equation}\label{eqw}
\left\|u\right\|^{p}_{L^{p}(B(x,r))}=\int\limits_{\left\{y\in I_{e}\,:\,\left|x-y\right|<r\right\}}\left|u_{e}(y)\right|^{p}dy+\sum^{N}_{\stackrel{j=1}{j\neq e}}\int\limits_{\left\{y\in I_{j}\,:\,x+y<r\right\}}\left|u_{j}(y)\right|^{p}dy.
\end{equation}

For  each minimizing sequence $\left\{(u_{n},v_{n})\right\}\subset X$ of $J(\alpha,\beta)$ we introduce  the following
sequence  (L\'evy concentration functions) $\rho_{n}:\left[0,\infty\right)\rightarrow\left[0,\alpha+\beta\right]$ by
	$$\rho_{n}(r)=\sup_{x\in\Gamma}\left(\left\|u_{n}\right\|^{2}_{B(x,r)}+\left\|v_{n}\right\|^{2}_{B(x,r)}\right).$$
Since $\left\{(u_{n},v_{n})\right\}$ is a minimizing sequence, then $\{\rho_{n}\}$ is a uniformly bounded sequence of nondecreasing functions on $\left[0,\infty\right)$. Moreover, it has a subsequence which is still denoted by  $\{\rho_{n}\}$ that converges pointwise and uniformly on compact sets to a nonnegative nondecreasing function $\rho(r):\left[0,\infty\right)\rightarrow\left[0,\alpha+\beta\right]$ (see  \cite[Lemma 3.2]{NQTG}). Define
\begin{equation*}
\tau=\lim_{r\rightarrow\infty}\rho(r)=\lim_{r\rightarrow\infty}\lim_{n\rightarrow\infty}\rho_{n}(r)=\lim_{r\rightarrow\infty}\lim\limits_{n\rightarrow\infty}\sup_{x\in\Gamma}\left(\left\|u_{n}\right\|^{2}_{B(x,r)}+\left\|v_{n}\right\|^{2}_{B(x,r)}\right).
\end{equation*}
Since $\left\|u_{n}\right\|^{2}_{2}=\alpha$ and $\left\|v_{n}\right\|^{2}_{2}=\beta$, it is clear that $\tau\in\left[0,\alpha+\beta\right]$. By concentration-compactness lemma for star graphs \cite[Lemma 3.3]{NQTG}, we have three (mutually exclusive) possibilities.\\
$(I)$\,(\textit{Compactness}) $\tau=\alpha+\beta$. Then, up to a subsequence, at least one of the two following cases occur:\\
\,$(I_{1})$(\textit{Convergence}) There exists $(\Phi_1,\Phi_2)\in X$ such that $u_{n}\rightarrow\Phi_1$ and $v_{n}\rightarrow\Phi_2$ in $L^{p}(\Gamma)$ as $n\rightarrow\infty$ for all $2\leq p\leq\infty$.\\
\,$(I_{2})$\,(\textit{Runaway}) There exists $j^{*}$ such that for any $e\neq j^*$, $r>0$, and $2\leq p\leq\infty$
\begin{equation*}
\|(u_{n})_{e}\|_{p}\rightarrow0, \quad\|(v_{n})_{e}\|_{p}\rightarrow0, \quad\left\|u_{n}\right\|_{L^{p}(B(0,r))}\rightarrow0, \quad\left\|v_{n}\right\|_{L^{p}(B(0,r))}\rightarrow0
\end{equation*}
as $n\to \infty$.

\noindent$(II)$\,(\textit{Vanishing}) $\tau=0$. Then, up to a subsequence, $u_{n}\rightarrow0$ and $v_{n}\rightarrow0$ in $L^{p}(\Gamma)$ as $n\rightarrow\infty$ for all $2< p\leq\infty$.\\
$(III)$\,(\textit{Dichotomy}) $0<\tau<\alpha+\beta$.

In what follows we show that case $(I_1)$ holds ruling out consequently $(II), (III)$, and $(I_2)$.
The following two propositions are used  to show that vanishing case does not hold. 
\begin{proposition}\label{lemma2}
Let $\left\{(u_{n},v_{n})\right\}\subset X$ be a minimizing sequence for $J(\alpha,\beta)$. Then there exist constants $A>0$ and $\lambda>0$ such that

$(i)$\, $\left\|u_{n}\right\|_{H^1}+\left\|v_{n}\right\|_{H^1}\leq A$ for all $n\in\mathbb{N}$;

$(ii)$\, $\left\|u_{n}\right\|^{4}_{4}+\left\|v_{n}\right\|^{4}_{4}\geq\lambda$ for $n$ large enough.
\end{proposition}
\begin{proof}
$(i)$ From  \eqref{Gagli1}-\eqref{Gagli2}, inequality $J(\alpha, \beta)<0$, and $\|u_n\|_2^2=\alpha, \|v_n\|_2^2=\beta$ we obtain
\begin{equation*}
\begin{split}
&\left\|u_{n}\right\|^{2}_{H^1}+\left\|v_{n}\right\|^{2}_{H^1}=E(u_{n},v_{n})+\gamma\left(\left|u_{n1}(0)\right|^{2}+\left|v_{n1}(0)\right|^{2}\right)+\frac{1}{2}\left(a\left\|u_{n}\right\|^{4}_{4}+c\left\|v_{n}\right\|^{4}_{4}\right)\\
&+b\|u_{n}v_{n}\|^{2}_2+\alpha+\beta
\leq\sup_{n}E(u_{n},v_{n})+C\left(\left\|u'_{n}\right\|_{2}+\left\|v'_{n}\right\|_{2}\right)+\alpha+\beta\\
&\leq C\left(1+\left\|u_{n}\right\|_{H^1}+\left\|v_{n}\right\|_{H^1}\right).
\end{split}
\end{equation*}
Since quadratic term $\left\|u_{n}\right\|^{2}_{H^1}+\left\|v_{n}\right\|^{2}_{H^1}$ is less than the linear one, the  existence of the desired bound $A$ follows.\\
To prove $(ii)$, let us assume that such $\lambda$ does not exist, then
\begin{equation}\label{qq1}
\liminf\limits_{n\rightarrow\infty}\left(\left\|u_{n}\right\|^{4}_{4}+\left\|v_{n}\right\|^{4}_{4}\right)=0.
\end{equation}
It follows from  Gagliardo-Nirenberg inequality \eqref{GNIq} (for $p=4, q=\infty$) and $(i)$ that $\left(\left|u_{n1}(0)\right|+\left|v_{n1}(0)\right|\right)\leq C\left(\left\|u_{n}\right\|^{2/3}_{4}+\left\|v_{n}\right\|^{2/3}_{4}\right)$. Hence, by \eqref{qq1}, we obtain $$\liminf\limits_{n\rightarrow\infty}\left(\left|u_{n1}(0)\right|+\left|v_{n1}(0)\right|\right)=0.$$ It is easily seen that   there exists $C>0$ such that 
\begin{equation*}
\begin{split}
0\leq E(u_{n},v_{n})+C\big(\left\|u_{n}\right\|^{4}_{4}+\left\|v_{n}\right\|^{4}_{4}+\gamma\left(\left|u_{n1}(0)\right|^{2}+\left|v_{n1}(0)\right|^{2}\right)\big),
\end{split}
\end{equation*}
then $J(\alpha,\beta)=\lim\limits\limits_{n\rightarrow\infty}E(u_{n},v_{n})\geq0$ which  contradicts  $J(\alpha, \beta)<0$. 
\end{proof}
\begin{proposition}\label{lemma3}
Let $u,v\in H^{1}(\Gamma)$, $A, \lambda>0$. Assume $\left\|u\right\|_{H^1}+\left\|v\right\|_{H^1}\leq A$ and $\left\|u\right\|^{4}_{4}+\left\|v\right\|^{4}_{4}\geq\lambda$, then there exists $B=B(A,\lambda)$ such that
\begin{equation*}
\sup_{x\in\Gamma}\left(\left\|u\right\|^{4}_{L^{4}(B(x,\frac{1}{2}))}+\left\|v\right\|^{4}_{L^{4}(B(x,\frac{1}{2}))}\right)\geq B.
\end{equation*}
\end{proposition}
\begin{proof}
Since $\left\|u\right\|^{4}_{4}+\left\|v\right\|^{4}_{4}\geq\lambda$, then we can choose $e_0\in\left\{1,\ldots,N\right\}$ and $\lambda_{e_0}>0$ such that $\|u_{e_0}\|^{4}_{4}+\|v_{e_0}\|^{4}_{4}\geq \lambda_{e_0}$. Let $n\in I_{e_0}$ be a natural number and  $f\in L^{p}(\Gamma)$, then 
\begin{equation*}
\left\|f\right\|^{p}_{L^{p}(B(n,\frac{1}{2}))}=\begin{cases} 
\int^{n+\frac{1}{2}}_{n-\frac{1}{2}}\left|f_{e_0}(x)\right|^{p}dx\,\,\,\mbox{for}\,\,\,n\geq1\\
\sum\limits_{e=1}^N\int^{\frac{1}{2}}_{0}\left|f_{e}(x)\right|^{p}dx\,\,\,\mbox{for}\,\,\,n=0.
\end{cases}
\end{equation*}
Hence
\begin{equation*}
\begin{split}
\|u_{e_0}\|^{2}_{H^{1}}+\|v_{e_0}\|^{2}_{H^{1}}&\leq\sum_{n\in\mathbb{N}_0}\left\{\left\|u\right\|^{2}_{H^1(B(n,\frac{1}{2}))}+\left\|v\right\|^{2}_{H^1(B(n,\frac{1}{2}))}\right\}\\
&\leq \frac{A^{2}}{\|u_{e_0}\|^{4}_{4}+\|v_{e_0}\|^{4}_{4}}\left(\|u_{e_0}\|^{4}_{4}+\|v_{e_0}\|^{4}_{4}\right)\\
&\leq\sum_{n\in\mathbb{N}_0}\frac{A^{2}}{\|u_{e_0}\|^{4}_{4}+\|v_{e_0}\|^{4}_{4}}\left(\left\|u\right\|^{4}_{L^{4}(B(n,\frac{1}{2}))}+\left\|v\right\|^{4}_{L^{4}(B(n,\frac{1}{2}))}\right).
\end{split}
\end{equation*}
Therefore, there must exist $n_{0}\in\mathbb{N}_0$ such that
\begin{equation}\label{left_estim}
\begin{split}
&\left\|u\right\|^{2}_{H^1(B(n_0,\frac{1}{2}))}+\left\|v\right\|^{2}_{H^1(B(n_0,\frac{1}{2}))}\\&\leq\frac{A^{2}}{\|u_{e_0}\|^{4}_{4}+\|v_{e_0}\|^{4}_{4}}\left(\left\|u\right\|^{4}_{L^{4}(B(n_{0},\frac{1}{2}))}+\left\|v\right\|^{4}_{L^{4}(B(n_{0},\frac{1}{2}))}\right).
\end{split}
\end{equation}
Now, from  inequality \eqref{GN_comp} we have
\begin{equation}\label{right_estim}
\begin{split}
&\left\|u\right\|^{4}_{L^{4}(B(n_{0},\frac{1}{2}))}+\left\|v\right\|^{4}_{L^{4}(B(n_{0},\frac{1}{2}))}\leq C (\|u\|_{H^1(B(n_0,\frac{1}{2}))}\|u\|_2^3+\|v\|_{H^1(B(n_0,\frac{1}{2}))}\|v\|_2^3)\\&\leq C(\|u\|_{H^1(B(n_0,\frac{1}{2}))}^2\|u\|_2^2+\|v\|_{H^1(B(n_0,\frac{1}{2}))}^2\|v\|_2^2) \leq C\left(\left\|u\right\|^{2}_{H^1(B(n_0,\frac{1}{2}))}+\left\|v\right\|^{2}_{H^1(B(n_0,\frac{1}{2}))}\right)^{2},
\end{split}
\end{equation}
where $C$ does not depend  on $u$ and $v$. Then, combining \eqref{left_estim} and \eqref{right_estim}, we get 
\begin{align*}
\left\|u\right\|^{2}_{H^1(B(n_0,\frac{1}{2}))}+\left\|v\right\|^{2}_{H^1(B(n_0,\frac{1}{2}))}\geq  \frac{\|u_{e_0}\|^{4}_{4}+\|v_{e_0}\|^{4}_{4}}{C A^{2}}.
\end{align*}

Thus, from \eqref{left_estim} we obtain
\begin{equation*}
B:=\frac{\lambda_{e_0}^{2}}{CA^{4}}\leq\left\|u\right\|^{4}_{L^{4}(B(n_{0},\frac{1}{2}))}+\left\|v\right\|^{4}_{L^{4}(B(n_{0},\frac{1}{2}))}.
\end{equation*}
\end{proof}

The following lemma rules out the vanishing case.

\begin{lemma} Let  $\left\{\left(u_{n},v_{n}\right)\right\}$  be a minimizing sequence for $J(\alpha,\beta)$, then  $\tau>0$.
\end{lemma}
\begin{proof}
It follows from Propositions \ref{lemma2} and \ref{lemma3} that there must exist a sequence $\left\{x_{n}\right\}\in \mathbb{R}^{+}$ and $B>0$ such that
\begin{equation*}
\left\|u_{n}\right\|^{4}_{L^{4}(B(x_{n},1/2))}+\left\|v_{n}\right\|^{4}_{L^{4}(B(x_{n},1/2))}\geq B
\end{equation*}
for all $n$. Hence, using \eqref{eqw} and the  Sobolev embedding of $H^{1}(\Gamma)$ into $L^{\infty}(\Gamma)$, we get
\begin{equation*}
\begin{split}
B&\leq \left\|u_{n}\right\|^{2}_{\infty}\left\|u_{n}\right\|^{2}_{L^{2}(B(x_{n},1/2))}+\left\|v_{n}\right\|^{2}_{\infty}\left\|v_{n}\right\|^{2}_{L^{2}(B(x_{n},1/2))}\\
&\leq CA^2\left(\left\|u_{n}\right\|^{2}_{L^{2}(B(x_{n},1/2))}+\left\|v_{n}\right\|^{2}_{L^{2}(B(x_{n},1/2))}\right),
\end{split}
\end{equation*}
where $A$ is from Proposition \ref{lemma2}.
Thus,
\begin{equation*}
\tau=\lim_{r\rightarrow\infty}\rho(r)\geq \rho (\frac{1}{2})=\lim\limits_{n\rightarrow\infty}\rho_{n}(\frac{1}{2})\geq\frac{B}{CA^2}>0.
\end{equation*}
\end{proof}
To exclude the possibility of dichotomy case, we first define
\begin{equation*}\label{E12}
\begin{split}
E_{1}(u)=\|u'\|_2^2-\gamma|u_{1}(0)|^2-\frac{b^{2}-ac}{2(b-c)}\left\|u\right\|^{4}_{4}, \qquad E_{2}(u)=\|u'\|_2^2-\gamma|u_{1}(0)|^2-\frac{b^{2}-ac}{2(b-a)}\left\|u\right\|^{4}_{4},
\end{split}
\end{equation*}
and  minimization problems:
\begin{equation}\label{J12}
\begin{split}
&J_{1}(\alpha)=\inf\left\{E_{1}(u):\,u\in H^{1}(\Gamma),\,\, \left\|u\right\|^{2}_{2}=\alpha \right\}, \\
&J_{2}(\beta)=\inf\left\{E_{2}(u):\,u\in H^{1}(\Gamma),\,\, \left\|u\right\|^{2}_{2}=\beta \right\}.
\end{split}
\end{equation}
The corresponding sets of minimizers for $J_{1}$ and $J_{2}$ are denoted by 
\begin{equation*}\label{G123}
\begin{split}
\mathcal{G}_{1}(\alpha)&=\left\{u\in H^{1}(\Gamma)\,:\,J_{1}(\alpha)=E_{1}(u),\, \left\|u\right\|^{2}_{2}=\alpha \right\},\\
\mathcal{G}_{2}(\beta)&=\left\{u\in H^{1}(\Gamma)\,:\,J_{2}(\alpha)=E_{2}(u),\, \left\|u\right\|^{2}_{2}=\beta \right\}.
\end{split}
\end{equation*}
The following four technical propositions will be used to rule out the possibility of dichotomy of minimizing sequences.
\begin{proposition}\label{lemma4}
Let $N\sqrt{\omega}-\gamma\leq\frac{2\gamma}{N}$, then the following equalities hold:

$(i)$\, $J_{1}(\alpha)=-\frac{1}{N^{2}}\left(\frac{1}{3}\left[\frac{b^{2}-ac}{2(b-c)}\right]^{2}\alpha^{3}+\frac{b^{2}-ac}{2(b-c)}\gamma\alpha^{2}+\gamma^{2}\alpha\right)$;

$(ii)$ $J_{2}(\beta)=-\frac{1}{N^{2}}\left(\frac{1}{3}\left[\frac{b^{2}-ac}{2(b-a)}\right]^{2}\beta^{3}+\frac{b^{2}-ac}{2(b-a)}\gamma\beta^{2}+\gamma^{2}\beta\right)$;

$(iii)$ $J(\alpha,\beta)=J_{1}(\alpha)+J_{2}(\beta)$.

\end{proposition}
\begin{proof}
$(i)\, $Since $\gamma>0$, $\omega>\frac{\gamma^{2}}{N^{2}}$, and $N\sqrt{\omega}-\gamma\leq\frac{2\gamma}{N}$, by  \cite[Theorem 2]{NQTG} (see also Remark \ref{rem_minimizer}), we have  
\begin{equation*}
\mathcal{G}_{1}(\alpha)=\left\{e^{i\theta_{1}}\sqrt{\frac{b-c}{b^{2}-ac}}\Phi_{\omega,\gamma}(x)\,:\,\theta_{1}\in\mathbb{R}\right\}
\end{equation*}
and
\begin{equation*}
\mathcal{G}_{2}(\beta)=\left\{e^{i\theta_{2}}\sqrt{\frac{b-a}{b^{2}-ac}}\Phi_{\omega,\gamma}(x)\,:\,\theta_{2}\in\mathbb{R}\right\},
\end{equation*}
where, for $e=1,\ldots,N$, $\left(\Phi_{\omega,\gamma}\right)_{e}=\phi_{\omega,\gamma}$, with
\begin{equation*}
\phi_{\omega,\gamma}(x)=\sqrt{2\omega}\,\mbox{sech}\left(\sqrt{\omega}x+\mbox{tanh}^{-1}\left(\frac{\gamma}{N\sqrt{\omega}}\right)\right).
\end{equation*} 
Therefore, $J_{1}(\alpha)=E_{1}(\sqrt{\frac{b-c}{b^{2}-ac}}\Phi_{\omega,\gamma})$ and $J_{2}(\beta)=E_{2}(\sqrt{\frac{b-a}{b^{2}-ac}}\Phi_{\omega,\gamma})$. Hence
\begin{equation}\label{J_1_alpha}
\begin{split}
&J_{1}(\alpha)=N\frac{b-c}{b^{2}-ac}\left(\int_{\mathbb{R}^{+}}\left|\phi'_{\omega,\gamma}(x)\right|^{2}dx-\frac{1}{2}\int_{\mathbb{R}^{+}}\left|\phi_{\omega,\gamma}(x)\right|^{4}dx-\frac{\gamma}{N}\left|\phi_{\omega,\gamma}(0)\right|^{2}\right)\\
&=2N\frac{b-c}{b^{2}-ac}\Bigg(\omega^{2}\int_{\mathbb{R}^{+}}\mbox{sech}^{2}\left(\sqrt{\omega}x+\mbox{tanh}^{-1}\left(\frac{\gamma}{N\sqrt{\omega}}\right)\right)dx\\
&-2\omega^{2}\int_{\mathbb{R}^{+}}\mbox{sech}^{4}\left(\sqrt{\omega}x+\mbox{tanh}^{-1}\left(\frac{\gamma}{N\sqrt{\omega}}\right)\right)dx-\frac{\omega\gamma}{N}+\frac{\gamma^{3}}{N^{3}}\Bigg)\\&=2N\frac{b-c}{b^{2}-ac}\Bigg(\frac{\omega^2}{\sqrt{\omega}}(1-\frac{\gamma}{N\sqrt{\omega}})-\frac{2\omega^2}{\sqrt{\omega}}(1-\frac{\gamma}{N\sqrt{\omega}})+\frac{2\omega^2}{3\sqrt{\omega}}(1-\frac{\gamma^3}{N^3\omega^{3/2}})-\frac{\omega\gamma}{N}+\frac{\gamma^3}{N^3}\Bigg)\\
&=-\frac{2}{3}\frac{b-c}{b^{2}-ac}\left(N\omega^{3/2}-\frac{\gamma^{3}}{N^{2}}\right).
\end{split}
\end{equation}
Recalling that $\alpha=2\frac{b-c}{b^2-ac}(N\sqrt{\omega}-\gamma)$, we get 
\begin{equation*}
\begin{split}
&J_1(\alpha)=-\frac{2}{3}\frac{b-c}{b^{2}-ac}\left(N\omega^{3/2}-\frac{\gamma^{3}}{N^{2}}\right)\\&=-\frac{1}{N^{2}}\left(\frac{2}{3}\frac{b-c}{b^{2}-ac}\Big[(N\sqrt{\omega}-\gamma)^3+3\Big(\gamma(N\sqrt{\omega}-\gamma)^2+\gamma^2(N\sqrt{\omega}-\gamma)\Big)\Big]\right)\\
&=-\frac{1}{N^{2}}\left(\frac{1}{3}\left[\frac{b^{2}-ac}{2(b-c)}\right]^{2}\alpha^{3}+\frac{b^{2}-ac}{2(b-c)}\gamma\alpha^{2}+\gamma^{2}\alpha\right).
\end{split}
\end{equation*}
In the same manner we can show $(ii)$.\\
$(iii)$\, Next we prove  $J(\alpha,\beta)\geq J_{1}(\alpha)+J_{2}(\beta)$. Using the Cauchy-Schwartz inequality and the Young inequality, we have
\begin{equation*}
\int_{\Gamma}\left|u\right|^{2}\left|v\right|^{2}dx\leq \sqrt{\frac{b-a}{b-c}}\left\|u\right\|^{2}_{4}\sqrt{\frac{b-c}{b-a}}\left\|v\right\|^{2}_{4}\leq\frac{1}{2}\left(\frac{b-a}{b-c}\left\|u\right\|^{4}_{4}+\frac{b-c}{b-a}\left\|u\right\|^{4}_{4}\right).
\end{equation*}
Hence
\begin{equation}\label{EEQQ1}
\begin{split}
&E(u,v)\geq\|u'\|_2^2+\|v'\|_2^2-\gamma\left(|u_{1}(0)|^2+|v_{1}(0)|^2\right)\\&-\frac{1}{2}\left(a\left\|u\right\|^{4}_{4}+c\left\|v\right\|^{4}_{4}\right)-\frac{b}{2}\left(\frac{b-a}{b-c}\left\|u\right\|^{4}_{4}+\frac{b-c}{b-a}\left\|v\right\|^{4}_{4}\right)\\
&=\|u'\|_2^2+\|v'\|_2^2-\gamma\left(|u_{1}(0)|^2+|v_{1}(0)|^2\right)-\frac{1}{2}\left(\frac{b^{2}-ac}{b-c}\left\|u\right\|^{4}_{4}+\frac{b^{2}-ac}{b-a}\left\|v\right\|^{4}_{4}\right)\\
&=E_{1}(u)+E_{2}(v).
\end{split}
\end{equation}
Taking the infimum on $u$ and $v$, we obtain $J(\alpha,\beta)\geq J_{1}(\alpha)+J_{2}(\beta)$. On the other hand, observe that $\left\|\sqrt{\frac{b-c}{b^{2}-ac}}\Phi_{\omega,\gamma}\right\|^{2}_{2}=\alpha$ and $\left\|\sqrt{\frac{b-a}{b^{2}-ac}}\Phi_{\omega,\gamma}\right\|^{2}_{2}=\beta$ for any $\omega>0$ fixed. Thus, we have
\begin{equation*}
\begin{split}
J(\alpha,\beta)&\leq E\left(\sqrt{\frac{b-c}{b^{2}-ac}}\Phi_{\omega,\gamma},\sqrt{\frac{b-a}{b^{2}-ac}}\Phi_{\omega,\gamma}\right)\\
&=E_{1}\left(\sqrt{\frac{b-c}{b^{2}-ac}}\Phi_{\omega,\gamma}\right)+ E_{2}\left(\sqrt{\frac{b-a}{b^{2}-ac}}\Phi_{\omega,\gamma}\right)\\
&=J_{1}(\alpha)+J_{2}(\beta)\leq J(\alpha,\beta).
\end{split}
\end{equation*}
Hence we arrive at $(iii)$, and the result is proved.
\end{proof}
\begin{proposition}\label{lemma5}
Let $\delta_{1} \in (0,\alpha)$ and $\delta_{2} \in (0,\beta)$, then 
\begin{equation*}
J_{1}(\alpha)<J_{1}(\delta_{1})+J_{1}(\alpha-\delta_{1})\,\,\,\,\,\,\mbox{and}\,\,\,\,\,\,J_{2}(\beta)<J_{2}(\delta_{2})+J_{2}(\beta-\delta_{2}).
\end{equation*}
\end{proposition}
\begin{proof}
We claim that if $\eta>1$ and $\alpha>0$, then $J_{1}(\eta\alpha)<\eta J_{1}(\alpha)$. To prove this inequality, consider  a minimizing sequence $\left\{u_{n}\right\}$ for $J_{1}(\alpha)$, and set $\widetilde{u}_{n}=\sqrt{\eta}u_{n}$ for all $n$. Then $\left\|\widetilde{u}_{n}\right\|^{2}_{2}=\eta\alpha$ and hence 
\begin{equation}\label{eqqq}
\begin{split}
&J_{1}(\eta\alpha)\leq E_{1}(\widetilde{u}_{n})=\eta \|u'_n\|_2^2-\gamma\eta|u_{n1}(0)|^2-\eta^{2}\frac{b^{2}-ac}{2(b-c)}\left\|u_{n}\right\|^{4}_{4}\\&=\eta E_{1}(u_{n})-\left(\eta^{2}-\eta\right)\frac{b^{2}-ac}{2(b-c)}\left\|u_{n}\right\|^{4}_{4}.
\end{split}
\end{equation}
On the other hand, repeating the proof of Proposition \ref{lemma2}, we can show that there exists a constant $\lambda>0$ such that $\left\|u_{n}\right\|^{4}_{4}\geq\lambda$ for $n$ large enough. Hence, taking $n\rightarrow\infty$ in \eqref{eqqq}, we obtain
\begin{equation*}
J_{1}(\eta\alpha)\leq \eta J_{1}(\alpha)-\left(\eta^{2}-\eta\right)\frac{b^{2}-ac}{2(b-c)}\lambda<\eta J_{1}(\alpha).
\end{equation*}
Without  loss of generality, we suppose that  $\delta_{1}>\alpha-\delta_{1}$.  Then, by the claim proved above, we have
\begin{equation*}
\begin{split}
J_{1}(\alpha)&=J_{1}(\delta_{1}(1+\frac{\alpha-\delta_{1}}{\delta_{1}}))<\left(1+\frac{\alpha-\delta_{1}}{\delta_{1}}\right)J_{1}(\delta_{1})\\
&<J_{1}(\delta_{1})+\frac{\alpha-\delta_{1}}{\delta_{1}}\left(\frac{\delta_{1}}{\alpha-\delta_{1}}J_{1}(\alpha-\delta_{1})\right)=J_{1}(\delta_{1})+J_{1}(\alpha-\delta_{1}).
\end{split}
\end{equation*}
Also, if $\delta_{1}=\alpha-\delta_{1}$, then we get
\begin{equation*}
J_{1}(\alpha)=J_{1}(2\delta_{1})<2J_{1}(\delta_{1})=J_{1}(\delta_{1})+J_{1}(\alpha-\delta_{1}).
\end{equation*}
The result for $J_{2}$ can be derived in an  analogous  way.
\end{proof}
The following proposition states strict subadditivity of $J$.
\begin{proposition}\label{lemma6}
Let $\delta_{1}\in\left[0,\alpha\right]$ and $\delta_{2}\in\left[0,\beta\right]$ be such that $0<\delta_{1}+\delta_{2}<\alpha+\beta$. Then
\begin{equation*}
J(\alpha,\beta)<J(\delta_{1},\delta_{2})+J(\alpha-\delta_{1},\beta-\delta_{2}).
\end{equation*}
\end{proposition}
\begin{proof}
We divide the proof into the following three cases:\\
\textbf{Case $1$.} Let $\delta_{1} \in (0,\alpha)$ and $\delta_{2} \in (0,\beta)$.  Using  Proposition \ref{lemma4}-$(iii)$ and Proposition \ref{lemma5}, we obtain
\begin{equation}\label{EQQ4}
J(\alpha,\beta)=J_{1}(\alpha)+J_{2}(\beta)<J_{1}(\delta_{1})+J_{1}(\alpha-\delta_{1})+J_{2}(\delta_{2})+J_{2}(\alpha-\delta_{2}).
\end{equation}
On the other hand, from \eqref{EEQQ1}  we get
\begin{equation}\label{EQQ5}
J_{1}(\delta_{1})+J_{2}(\delta_{2})\leq J(\delta_{1},\delta_{2})\,\,\,\,\,\,\mbox{and}\,\,\,\,\,\,J_{1}(\alpha-\delta_{1})+J_{2}(\beta-\delta_{2})\leq J(\alpha-\delta_{1},\beta-\delta_{2}).
\end{equation}
Thus, combining \eqref{EQQ4} and \eqref{EQQ5}, we obtain
\begin{equation*}
J(\alpha,\beta)<J(\delta_{1},\delta_{2})+J(\alpha-\delta_{1},\beta-\delta_{2}),
\end{equation*}
as desired.\\
\textbf{Case $2$.} Assume that $\delta_{1}=0$ and $\delta_{2}\in\left(0,\beta\right]$.  We consider the following variational problem
\begin{equation}\label{mini}
J(0,\delta_{2})=\inf\left\{\|v'\|_2^2-\gamma|v_{1}(0)|^2-\frac{c}{2}\left\|v\right\|^{4}_{4}\,:\,v\in H^{1}(\Gamma),\,\left\|v\right\|^{2}_{2}=\delta_{2}\right\}.
\end{equation}
For $c>0$ the minimizer of \eqref{mini} is given by $\Phi_2(x)=\left(\phi_2(x)\right)_{e=1}^N$, where
\begin{equation*}
\phi_2(x)=\sqrt{\frac{2\omega_{2}}{c}}\mbox{sech}\,\left(\sqrt{\omega_{2}}x+\mbox{tanh}^{-1}\left(\frac{\gamma}{N\sqrt{\omega_{2}}}\right)\,\right)
\end{equation*}
and $\omega_{2}=\left(\frac{c\delta_{2}+2\gamma}{2N}\right)^{2}$ (see formulas (5.4) and (5.1) in \cite{NQTG}). Therefore,
\begin{equation*}\label{QQ0}
\begin{split}
J(0,\delta_{2})=N\|\phi'_2\|^{2}_{2}-\gamma\left|\phi_2(0)\right|^{2}-N\frac{c}{2}\|\phi_2\|^{4}_{4}=-\frac{1}{N^{2}}\left(\frac{c^{2}}{12}\delta^{3}_{2}+\frac{c}{2}\gamma\delta^{2}_{2}+\gamma^{2}\delta_{2}\right).
\end{split}
\end{equation*}
Now, from  Proposition \ref{lemma4}-$(ii)$, we have
\begin{equation*}
J_{2}(\delta_{2})=-\frac{1}{N^{2}}\left(\frac{1}{3}\left[\frac{b^{2}-ac}{2(b-a)}\right]^{2}\delta^{3}_{2}+\frac{b^{2}-ac}{2(b-a)}\gamma\delta^{2}_{2}+\gamma^{2}\delta_{2}\right).
\end{equation*}
Hence $J(0,\delta_{2})>J_{2}(\delta_{2})$. Then
\begin{equation*}
J(\alpha,\beta)=J_{1}(\alpha)+J_{2}(\beta)\leq J_{1}(\alpha)+J_{2}(\delta_{2})+J_{2}(\beta-\delta_{2})<J(0,\delta_{2})+J(\alpha,\beta-\delta_{2}).
\end{equation*}
\iffalse
Assume that $c<0$, then $$J(0,\delta_2)\geq \inf\{\|v'\|_2^2-\gamma|v_{1}(0)|^2: v\in H^1(\Gamma), \|v\|_2^2=\delta_2\}\geq \inf\sigma(-\Delta_\Gamma)\delta_2= -\frac{\gamma^2}{N^2}\delta_2.$$
Obviously \fi
\textbf{Case $3$.} Assume that $\delta_{1}\in\left(0,\alpha\right]$ and $\delta_{2}=0$. Analogously to Case 2, we consider the following variational problem
\begin{equation*}\label{mini2}
J(\delta_{1},0)=\inf\left\{\|u'\|_2^2-\gamma|u_{1}(0)|^2-\frac{a}{2}\left\|u\right\|^{4}_{4}\,:\,u\in H^{1}(\Gamma),\,\left\|u\right\|^{2}_{2}=\delta_{1}\right\}.
\end{equation*}
Since $a>0$, the minimizer of \eqref{mini} is given by $\Phi_1(x)=\left(\phi_1(x)\right)_{e=1}^N$, where
\begin{equation*}
\phi_1(x)=\sqrt{\frac{2\omega_{1}}{a}}\mbox{sech}\,\left(\sqrt{\omega_{1}}x+\mbox{tanh}^{-1}\left(\frac{\gamma}{N\sqrt{\omega_{1}}}\right)\,\right),
\end{equation*}
where $\omega_{1}=\left(\frac{a\delta_{1}+2\gamma}{2N}\right)^{2}$. Therefore,
\begin{equation}\label{QQ00}
J(\delta_{1},0)=-\frac{1}{N^{2}}\left(\frac{a^{2}}{12}\delta^{3}_{1}+\frac{a}{2}\gamma\delta^{2}_{1}+\gamma^{2}\delta_{1}\right).
\end{equation}
Now, from  Proposition \ref{lemma4}-$(i)$ we have
\begin{equation*}
J_{1}(\delta_{1})=-\frac{1}{N^{2}}\left(\frac{1}{3}\left[\frac{b^{2}-ac}{2(b-c)}\right]^{2}\delta^{3}_{1}+\frac{b^{2}-ac}{2(b-c)}\gamma\delta^{2}_{1}+\gamma^{2}\delta_{1}\right).
\end{equation*}
Hence from \eqref{QQ00} we get that $J(\delta_{1},0)>J_{1}(\delta_{1})$. Then
\begin{equation*}
J(\alpha,\beta)=J_{1}(\alpha)+J_{2}(\beta)\leq J_{1}(\delta_{1})+ J_{1}(\alpha-\delta_{1})+J_{2}(\beta)<J(\delta_{1},0)+J(\alpha-\delta_{1},\beta).
\end{equation*}
\end{proof}
We introduce two sequences $\{(f_n, g_n)\}$ and $\{(h_n, l_n)\}$ associated with an arbitrary minimizing sequence $\{(u_n, v_n)\}$ in the following way.  

Let $\sigma, \kappa\in C^{\infty}_{0}(\mathbb{R}^+)$. We assume that $\sigma$ is supported on $[0,2]$, $\sigma\equiv1$ on $\left[0,1\right]$ and  $\sigma^{2}+\kappa^{2}=1$ on $\mathbb{R}^{+}$. Set $\sigma_{r}(x)=\sigma(\frac{x}{r})$ and $\kappa_{r}(x)=\kappa(\frac{x}{r})$ for $r>0$. Let $\epsilon>0$, for sufficiently large $r$ we have $\tau-\epsilon<\rho(r)\leq\rho(2r)\leq\tau$. We can choose $N$ large enough so that
\begin{equation*}
\tau-\epsilon<\rho_{n}(r)\leq\rho_{n}(2r)<\tau+\epsilon
\end{equation*}
 for all $n\geq N$. Consequently, for each $n\geq N$ we can find $x_{n}\in\Gamma$ such that
\begin{equation}\label{56}
\begin{split}
\left(\left\|u_{n}\right\|^{2}_{L^{2}(B(x_{n},r))}+\left\|v_{n}\right\|^{2}_{L^{2}(B(x_{n},r))}\right)>\tau-\epsilon,\\
\left(\left\|u_{n}\right\|^{2}_{L^{2}(B(x_{n},2r))}+\left\|v_{n}\right\|^{2}_{L^{2}(B(x_{n},2r))}\right)<\tau+\epsilon.
\end{split}
\end{equation}
Let $x\in\Gamma$, $d(x,x_{n})$ denotes the distance between $x$ and $x_{n}$ given by  \eqref{55}. Define $(f_{n},g_{n})$ and $(h_{n},l_{n})$ such that $f_{n,e}(x)=\sigma_{r}(d(x,x_{n}))u_{n,e}(x)$, $g_{n,e}(x)=\sigma_{r}(d(x,x_{n}))v_{n,e}(x)$, $h_{n,e}(x)=\kappa_{r}(d(x,x_{n}))u_{n,e}(x)$ and $l_{n,e}(x)=\kappa_{r}(d(x,x_{n}))v_{n,e}(x)$ for $e=1,\ldots,N$. Observe that $Q(f_{n})\geq\left\|u_{n}\right\|^{2}_{L^{2}(B(x_{n},r))}$ and $Q(g_{n})\geq\left\|v_{n}\right\|^{2}_{L^{2}(B(x_{n},r))}$. Since $\sigma\leq1$ and $\supp \sigma_r\subseteq\{x\in \Gamma: d(x, x_n)\leq 2r\}$, we have  $Q(f_{n})\leq\left\|u_{n}\right\|^{2}_{L^{2}(B(x_{n},2r))}$ and $Q(g_{n})\leq\left\|v_{n}\right\|^{2}_{L^{2}(B(x_{n},2r))}$.

Observe that $\alpha-\left\|u_{n}\right\|^{2}_{L^{2}(B(x_{n},2r))}\leq Q(h_{n})$ and $\beta-\left\|v_{n}\right\|^{2}_{L^{2}(B(x_{n},2r))}\leq Q(l_{n})$. Since $\kappa\leq1$  and $\{x\in \Gamma: d(x,x_n)>2r\}\subseteq\supp\kappa_r$,  we have  $Q(h_{n})\leq\alpha-\left\|u_{n}\right\|^{2}_{L^{2}(B(x_{n},r))}$ and $Q(l_{n})\leq\beta-\left\|v_{n}\right\|^{2}_{L^{2}(B(x_{n},r))}$. 

\begin{proposition}\label{lemma8} Let $\{(u_n, v_n)\}$ be a minimizing sequence for \eqref{JOTA} and sequences $\{(f_n, g_n)\}$, $\{(h_n, l_n)\}$ be defined above. Then for every $\epsilon>0$ there exists  $N>0$ such that   for  $n\geq N$:\\
$(i)$\, $\left|Q(f_{n})+Q(g_{n})-\tau\right|<\epsilon$,\\
$(ii)$\, $\left|Q(h_{n})+Q(l_{n})-\left(\alpha+\beta-\tau\right)\right|<\epsilon$,\\
$(iii)$\, $E(f_{n},g_{n})+E(h_{n},l_{n})\leq E(u_{n},v_{n})+C\epsilon$ 
for some $C>0$ independent of $n$.
\end{proposition}
\begin{proof}
The proof of $(i), (ii)$ is a direct consequence of \eqref{56}.
To prove $(iii)$ notice that for $\frac{1}{r}<\varepsilon$
\begin{equation*}
\begin{split}
&E(f_{n},g_{n})=\int_{\Gamma}\sigma^{2}_{r}\left[(u'_{n})^{2}+(v'_{n})^{2}-\frac{a}{2}\left|u_{n}\right|^{4}-\frac{c}{2}\left|u_{n}\right|^{4}-b\left|u_{n}\right|^{2}\left|v_{n}\right|^{2}\right]dx\\
&-\gamma\left|\sigma_{r}(x_{n})\right|^{2}\left(\left|u_{n1}(0)\right|^{2}+\left|v_{n1}(0)\right|^{2}\right)+\int_{\Gamma}\left(\sigma^{2}_{r}-\sigma^{4}_{r}\right)\left[\frac{a}{2}\left|u_{n}\right|^{4}+\frac{c}{2}\left|v_{n}\right|^{4}+b\left|u_{n}\right|^{2}\left|v_{n}\right|^{2}\right]dx\\
&+\int_{\Gamma}(\sigma'_{r})^{2}\Big[\left|u_{n}\right|^{2}+\left|v_{n}\right|^{2}+2\sigma'_{r}\sigma_{r}\left(u'_{n}u_{n}+v'_{n}v_{n}\right)\Big]dx\\
&\leq\int_{\Gamma}\sigma^{2}_{r}\left[(u'_{n})^{2}+(v'_{n})^{2}-\frac{a}{2}\left|u_{n}\right|^{4}-\frac{c}{2}\left|u_{n}\right|^{4}-b\left|u_{n}\right|^{2}\left|v_{n}\right|^{2}\right]dx\\
&-\gamma\left|\sigma_{r}(x_{n})\right|^{2}\left(\left|u_{n1}(0)\right|^{2}+\left|v_{n1}(0)\right|^{2}\right)+C\epsilon.
\end{split}
\end{equation*}
Indeed, observe that 
$\left|\sigma'_{r}\right|_{\infty}=\left|\sigma'\right|_{\infty}/r\leq C\epsilon$ (since $1/r<\varepsilon$). Under the integral $\sigma_r$ has to be read as  $(\sigma_r)_{e=1}^N$. Moreover, introducing 
\begin{equation*}
\left(\chi_{(B(x_{n},2r)\setminus B(x_{n},r))}\right)_{e}(x)=
\begin{cases} 
1\,\,\,\mbox{if}\,\,r\leq d(x,x_{n})<2r\\
0\,\,\,\mbox{otherwise},\\
\end{cases}
\end{equation*} 
and using \eqref{56}, we get
\begin{equation*}
\begin{split}
\int_{\Gamma}\left(\sigma^{2}_{r}-\sigma^{4}_{r}\right)\left|u_{n}\right|^{4}dx&\leq \|u_n\|_\infty^2\sum^{N}_{e=1}\int_{I_{e}}\left(\chi_{(B(x_{n},2r)\setminus B(x_{n},r))}\right)_{e}(x)\left|u_{n,e}(x)\right|^{2}dx\leq C\epsilon,
\end{split}
\end{equation*}
and
\begin{equation*}
\begin{split}
&\int_{\Gamma}\left(\sigma^{2}_{r}-\sigma^{4}_{r}\right)\left|u_{n}\right|^{2}\left|v_{n}\right|^{2}dx\\&\leq \|u_n\|_\infty^2 \sum^{N}_{e=1}\int_{I_{e}}\left(\chi_{(B(x_{n},2r)\setminus B(x_{n},r))}\right)_{e}(x)\left|v_{n,e}(x)\right|^{2}dx\leq C\epsilon,
\end{split}
\end{equation*}
where $C$ denotes  constant independent of $r$ and $n$. Similarly, we obtain
\begin{equation*}
\begin{split}
E(h_{n},l_{n})&\leq\int_{\Gamma}\kappa^{2}_{r}\left[(u'_{n})^{2}+(v'_{n})^{2}-\frac{a}{2}\left|u_{n}\right|^{4}-\frac{c}{2}\left|u_{n}\right|^{4}-b\left|u_{n}\right|^{2}\left|v_{n}\right|^{2}\right]dx\\
&-\gamma\left|\kappa_{r}(x_{n})\right|^{2}\left(\left|u_{n1}(0)\right|^{2}+\left|v_{n1}(0)\right|^{2}\right)+C\epsilon.
\end{split}
\end{equation*}
Thus, since $\sigma^{2}+\kappa^{2}\equiv1$ on $\mathbb{R}^{+}$, we get $(iii)$.
\end{proof}
Below we rule out the dichotomy of the  minimizing sequences.
\begin{lemma}
Let $\left\{\left(u_{n},v_{n}\right)\right\}$ be a  minimizing sequence for $J(\alpha,\beta)$. Then the case  of dichotomy cannot occur.
\end{lemma}
\begin{proof}
Assume that $\tau\in (0, \alpha+\beta)$.
Let $\{(f_n, g_n)\}$ and $\{(h_n, l_n)\}$ be the sequences from Proposition \ref{lemma8}, then, up to subsequences (using the fact that minimizing sequence is $L^2$-bounded), we have   
\begin{equation*}
\begin{split}
|Q(f_n)-\alpha+\delta_1|<\frac{1}{2n},& \,\, |Q(g_n)-\beta+\delta_2|<\frac{1}{2n}, \,\, |Q(h_n)-\delta_1|<\frac{1}{2n},\,\, |Q(l_n)-\delta_2|<\frac{1}{2n}, \\
& E(f_n, g_n)+E(h_n, l_n)\leq E(u_n, v_n)+\frac{1}{n},
\end{split}
    \end{equation*}
  where  $\delta_{1}\in \left[0,\alpha\right]$ and $\delta_{2}\in\left[0,\beta\right]$ and, by Proposition \ref{lemma8}, $\tau=\alpha+\beta-(\delta_1+\delta_2)$.
  For $\delta_1>0, \delta_2>0, \alpha-\delta_1>0, \beta-\delta_2>0$ we define 
\begin{equation*}
a_{1n}=\sqrt{\frac{\alpha-\delta_{1}}{Q(f_{n})}},\,\,\,\,a_{2n}=\sqrt{\frac{\beta-\delta_{2}}{Q(g_{n})}},\,\,\,\,b_{1n}=\sqrt{\frac{\delta_{1}}{Q(h_{n})}},\,\,\,\,b_{2n}=\sqrt{\frac{\delta_{2}}{Q(l_{n})}}.
\end{equation*}
 Thus, we have 
 \begin{align*}
   & a_{1n},a_{2n},b_{1n},b_{2n}\underset{n\to\infty}{\rightarrow} 1,\\
  &  Q(a_{1n}f_{n})=\alpha-\delta_{1},\,\,\, Q(a_{2n}g_{n})=\beta-\delta_{2},\,\,\,Q(b_{1n}h_{n})=\delta_{1},\,\,\,Q(b_{2n}l_{n})=\delta_{2}. \end{align*} Consequently
\begin{equation}\label{58}
\begin{split}
&J(\alpha-\delta_{1},\beta-\delta_{2})\leq E(a_{1n}f_{n},a_{2n}g_{n})\leq E(f_n, g_n)+o(1),\\
&J(\delta_{1},\delta_{2})\leq E(b_{1n}h_{n},b_{2n}l_{n})\leq E(h_n, l_n)+o(1).
\end{split}
\end{equation}
By  Proposition \ref{lemma8}-$(iii)$ and \eqref{58},  we get 
 $$J(\alpha-\delta_{1},\beta-\delta_{2})+J(\delta_{1},\delta_{2})\leq E(u_n, v_n)+o(1),$$
 which implies 
 $$J(\alpha-\delta_{1},\beta-\delta_{2})+J(\delta_{1},\delta_{2})\leq J(\alpha, \beta).$$
 This contradicts Proposition \ref{lemma6}.
 
 Now assume that one of the numbers $\delta_1, \delta_2, \alpha-\delta_1, \beta-\delta_2$ is zero. Without loss of generality, we can assume  $\delta_1=0$, then  $Q(h_n)\to 0$ and $\delta_2>0$ due to $\tau<\alpha+\beta$.
 Therefore, by the Gagliardo-Nirenberg inequality, 
 \begin{equation*}
\begin{split}
&E(h_n, l_n)+o(1)\geq\left(\left\|h'_{n}\right\|^{2}_{2}+\left\|l'_{n}\right\|^{2}_{2}-\gamma\left|l_{n1}(0)\right|^{2}-\frac{c}{2}\left\|l_{n}\right\|^{4}_{4}\right)+o(1)\\
&\geq\left(\left\|l'_{n}\right\|^{2}_{2}-\gamma\left|l_{n1}(0)\right|^{2}-\frac{c}{2}\left\|l_{n}\right\|^{4}_{4}\right)+o(1)\geq J(0,\delta_{2})+o(1).
\end{split}
\end{equation*}
 Combining this inequality with the first one from \eqref{58}, we arrive at the contradiction again.

\end{proof}
Since we eliminated  the vanishing and the dichotomy cases, it follows from the concentration-compactness lemma \cite{NQTG} that $\tau=\alpha+\beta$, that is, we have compactness case. It only remains to prove that the minimizing sequence is not \textit{runaway}.

\begin{lemma}\label{lemma10}
 Let $\left\{\left(u_{n},v_{n}\right)\right\}$ be a minimizing sequence for $J(\alpha,\beta)$. Then, up to subsequence, it converges in $X$ to some $\left(\Phi_1,\Phi_2\right)$, which is a minimizer for $J(\alpha,\beta)$, that is, $\left\|\Phi_1\right\|^{2}_{2}=\alpha$, $\left\|\Phi_2\right\|^{2}_{2}=\beta$, and $E(\Phi_1,\Phi_2)=J(\alpha,\beta)$.
\end{lemma}
\begin{proof} Set
\begin{equation*}
\begin{split}
&E_0(u,v)=\left\|u'\right\|^{2}_{2}+\left\|v'\right\|^{2}_{2}-\frac{1}{2}\left(a\left\|u\right\|^{4}_{4}+c\left\|v\right\|^{4}_{4}\right)-b\|uv\|^{2}_2,\\
&J_0(\alpha,\beta)=\inf\left\{E_0(u,v)\,:\,u,v\in H^{1}(\Gamma),\,\,\left\|u\right\|^{2}_{2}=\alpha,\,\,\left\|v\right\|^{2}_{2}=\beta\right\}.
\end{split}
\end{equation*}
By absurd, suppose that $\left\{(u_{n},v_{n})\right\}$ is the  \textit{runaway} sequence. Then from Proposition \ref{lemma2}-$(i)$ and the Gagliardo-Nirenberg inequality,   we have that $\left|u_{n,e}(0)\right|\rightarrow0$ and $\left|v_{n,e}(0)\right|\rightarrow0$, for $e \neq j^*$ ($j^*$ is from  definition $(I_2)$ of runaway sequence)  which implies that $\lim\limits_{n\rightarrow\infty}\left(E(u_{n},v_{n})-E_0(u_{n},v_{n})\right)=0$. Hence,
\begin{equation}\label{01}
J_0(\alpha,\beta)\leq J(\alpha,\beta).
\end{equation}
\iffalse
On the other hand, let $u,v\in H^{1}(\Gamma)$. Then, using the symmetric rearrangement theory for graphs, the following estimates hold (see Proposition \ref{Pr0p0} in Appendix)
\begin{equation*}
\left\|u\right\|_{2}=\left\|u^{*}\right\|_{2},\,\,\,\left\|u\right\|_{2}=\left\|u^{*}\right\|_{4},\,\,\,\left\|u'\right\|^{2}_{2}\geq\frac{4}{N^{2}}\left\|(u^{*})'\right\|^{2}_{2}
\end{equation*}
and
\begin{equation*}
\int_{\Gamma}\left|u\right|^{2}\left|v\right|^{2}dx\leq\int_{\Gamma}\left|u^{*}\right|^{2}\left|v^{*}\right|^{2}dx,
\end{equation*}
where $u^{*}$ and $v^{*}$ denote the symmetric rearrangements of $u$ and $v$ respectively.
\fi 
Take  $u,v\in H^{1}(\Gamma)$ such that $\left\|u\right\|^{2}_{2}=\alpha$ and $\left\|v\right\|^{2}_{2}=\beta$, and let  $u^{*},v^{*}\in H^{1}(\Gamma)$ be their symmetric rearrangements. Then, by  Proposition \ref{Pr0p0} in Appendix,  we obtain  $\left\|u^{*}\right\|^{2}_{2}=\alpha$, $\left\|v^{*}\right\|^{2}_{2}=\beta$, and $E_0(u,v)\geq E_0^*(u^*,v^*)$, where
\begin{equation*}
E_0^*(u,v)=\frac{4}{N^{2}}\left(\left\|(u^{*})'\right\|^{2}_{2}+\left\|(v^{*})'\right\|^{2}_{2}\right)-\frac{1}{2}\left(a\left\|u^{*}\right\|^{4}_{4}+c\left\|v^{*}\right\|^{4}_{4}\right)-b\|u^{*}v^{*}\|^2_2.
\end{equation*}
Since rearrangements maintain the mass constraint, the last inequality implies
\begin{equation}\label{symm}
J_0(\alpha,\beta)\geq\inf\left\{E_0^*(u,v):\,u,v\in H^{1}_{\mathrm{s}}(\Gamma),\,\,\left\|u\right\|^{2}_{2}=\alpha,\,\,\left\|v\right\|^{2}_{2}=\beta\right\},
\end{equation}
where $H^{1}_{\mathrm{s}}(\Gamma)=\{u\in H^1(\Gamma): u_1(x)=\ldots=u_N(x), x>0\}.$
It is easily seen that right-hand side of \eqref{symm} reduces to $N$ copies of the following problem on $\mathbb{R}^{+}$
\begin{equation*}
\inf\left\{E^*_{\mathbb{R}^+}(\psi,\varphi):\,\psi,\varphi\in H^{1}(\mathbb{R}^{+}),\,\,\|\psi\|^{2}_{2}=\frac{\alpha}{N},\,\,\|\varphi\|^{2}_{2}=\frac{\beta}{N}\right\},
\end{equation*}
where
\begin{equation*}
E^*_{\mathbb{R}^+}(\psi,\varphi)=\frac{4}{N^{2}}\left(\|\psi'\|^{2}_{2}+\|\varphi'\|^{2}_{2}\right)-\frac{1}{2}\left(a\|\psi\|^{4}_{4}+c\|\varphi\|^{4}_{4}\right)-b\|\psi\varphi\|^{2}_2.
\end{equation*}
Set for $\lambda > 0$ the rescaling  $\psi_{\lambda}(x)=\sqrt{\lambda}\psi(\lambda x)$ and $\varphi_{\lambda}(x)=\sqrt{\lambda}\varphi(\lambda x)$, then 
\begin{equation*}
E^*_{\mathbb{R}^+}(\psi_\lambda,\varphi_\lambda)=\frac{4\lambda^{2}}{N^{2}}\left(\|\psi'\|^{2}_{2}+\|\varphi'\|^{2}_{2}\right)-\frac{\lambda}{2}\left(a\|\psi\|^{4}_{4}+c\|\varphi\|^{4}_{4}\right)-b\lambda\|\psi\varphi\|^{2}_2.
\end{equation*}
Choosing $\lambda=\frac{N^{2}}{4}$, we get
\begin{equation*}
J_0(\alpha,\beta)\geq\frac{N^{3}}{4}d_{\mathbb{R}^+},
\end{equation*}
where
\begin{equation}\label{d_+}
d_{\mathbb{R}^+}=\inf\left\{E_{\mathbb{R}^{+}}(\psi,\varphi):\,\psi,\varphi\in H^{1}(\mathbb{R}^{+}),\,\,\|\psi\|^{2}_{2}=\frac{\alpha}{N},\,\,\|\varphi\|^{2}_{2}=\frac{\beta}{N}\right\},    
\end{equation}
\begin{equation*}
E_{\mathbb{R}^{+}}(\psi,\varphi)=\|\psi'\|^{2}_{2}+\|\varphi'\|^{2}_{2}-\frac{1}{2}\left(a\|\psi\|^{4}_{4}+c\|\varphi\|^{4}_{4}\right)dx-b\|\psi\varphi\|^{2}_2.
\end{equation*}
By  \cite[Theorem 2.1]{NguWa}, the solution to minimization problem \eqref{d_+} on $\mathbb{R}^{+}$ is  $\left(\psi_{\widetilde{\omega}},\varphi_{\widetilde{\omega}}\right)=\left(\sqrt{\frac{b-c}{b^{2}-ac}}\phi_{\widetilde{\omega}},\sqrt{\frac{b-a}{b^{2}-ac}}\phi_{\widetilde{\omega}}\right)$, where $\phi_{\widetilde{\omega}}(x)=\sqrt{2\widetilde{\omega}}\,\mbox{sech}(\sqrt{\widetilde{\omega}}x)$ and $\widetilde{\omega}$ is such that $\|\psi_{\widetilde{\omega}}\|^{2}_{2}=\frac{\alpha}{N}$ and $\|\varphi_{\widetilde{\omega}}\|^{2}_{2}=\frac{\beta}{N}$. Then  we get,
\begin{equation}\label{102}
J_0(\alpha,\beta)\geq\frac{N^{3}}{4} E_{\mathbb{R}^{+}}\left(\psi_{\widetilde{\omega}},\varphi_{\widetilde{\omega}}\right)=-\frac{2}{3}\left(\frac{N^{3}}{4}\widetilde{\omega}^{3/2}\right)\frac{2b-a-c}{b^{2}-ac}.
\end{equation}
By formula (3.3) in \cite{NguWa}, we have $\frac{\alpha}{N}=2\sqrt{\tilde{\omega}}\frac{b-c}{b^2-ac}$. Recalling that $\alpha=2\frac{b-c}{b^2-ac}(N\sqrt{\omega}-\gamma),$ we obtain $N\sqrt{\tilde{\omega}}=N\sqrt{\omega}-\gamma.$ Then from \eqref{102} we get
$$J_0(\alpha,\beta)\geq-\frac{2}{3}\left(\frac{N^{2}}{4}\widetilde{\omega}\left(N\sqrt{\omega}-\gamma\right)\right)\frac{2b-a-c}{b^{2}-ac}.$$

By Proposition \ref{lemma4}-$(iii)$ and formula \eqref{J_1_alpha}, we conclude
\begin{equation}\label{103}
\begin{split}
&J(\alpha,\beta)=-\frac{2}{3}\left(N\omega^{3/2}-\frac{\gamma^{3}}{N^{2}}\right)\frac{2b-a-c}{b^{2}-ac}\\&=-\frac{2}{3}\left(\omega\left(N\sqrt{\omega}-\gamma\right)+\gamma\left(\omega-\frac{\gamma^{2}}{N^{2}}\right)\right)\frac{2b-a-c}{b^{2}-ac}.
\end{split}
\end{equation}
Since $N\sqrt{\widetilde{\omega}}=(N\sqrt{\omega}-\gamma)\leq\frac{2\gamma}{N}$, then $\frac{N^{2}}{4}\widetilde{\omega}\leq\frac{\gamma^{2}}{N^{2}}<\omega$. Hence, by \eqref{102} and $\eqref{103}$, we deduce that $J_0(\alpha,\beta)>J(\alpha,\beta)$ which contradicts \eqref{01}. Therefore, $\left\{\left(u_{n},v_{n}\right)\right\}$ is not \textit{runaway}, and it  converges, up to subsequence, to  $(\Phi_1,\Phi_2)$ in  $L^{p}(\Gamma)\times L^{p}(\Gamma)$ for $p\geq2$ and weakly in $H^1(\Gamma)\times H^1(\Gamma)$. In particular, $Q(\Phi_1)=\alpha$ and $Q(\Phi_2)=\beta$, and, by the weak lower semicontinuity of the $H^{1}$ norm, we have 
\begin{equation*}
E(\Phi_1,\Phi_2)\leq \lim\limits_{n\rightarrow\infty}E(u_{n},v_{n})=J(\alpha,\beta),
\end{equation*}
whence $E(\Phi_1,\Phi_2)=J(\alpha,\beta)$ and $(\Phi_1,\Phi_2)\in \mathcal{G}(\alpha,\beta)$. Using  $E(\Phi_1,\Phi_2)= \lim\limits_{n\rightarrow\infty}E(u_{n},v_{n})$  and $\|u_n-\Phi_1\|_p\rightarrow 0$,  $\|v_n-\Phi_2\|_p\rightarrow 0$, $p\geq 2$, and the weak convergence,  we conclude 
$$\|(\Phi_1, \Phi_2)\|_X=\lim\limits_{n\rightarrow\infty}\|(u_{n},v_{n})\|_X.$$
Indeed, it is sufficient to observe that  
$$\|(u_{n},v_{n})\|_X^2=E(u_n, v_n)+\frac{1}{2}\left(a\|u_n\|_4^4+b\|v_n\|_4^4\right)+b\|u_nv_n\|^2_2+\gamma|u_{n1}(0)|^2+\gamma|v_{n1}(0)|^2.$$
Finally, since $X$ is the Hilbert space, we conclude that $(u_{n},v_{n})\rightarrow (\Phi_1,\Phi_2)$ in $X$. 
\end{proof}
\begin{remark}
In \cite{NguWa} the authors consider    minimizing problem on the line 
\begin{equation*}d_{\mathbb{R}}=\inf\{E_{\mathbb{R}}: (\psi, \varphi)\in H^1(\mathbb{R}), \,\|\psi\|_2^2=\tilde{\alpha},\, \|\varphi\|_2^2=\tilde{\beta}\},
\end{equation*}
where $E_{\mathbb{R}}=E_{\mathbb{R}^+}$, and  $\tilde{\alpha}=4 \sqrt{\tilde{\omega}} \frac{b-c}{b^{2}-a c}$, $\tilde{\beta}=4 \sqrt{\tilde{\omega}} \frac{b-a}{b^{2}-a c}$. Assuming $\tilde{\alpha}=\frac{2\alpha}{N}, \tilde{\beta}=\frac{2\beta}{N}$,  recalling that $d_{\mathbb{R}}=E_{\mathbb{R}}\left(\sqrt{\frac{b-c}{b^2-ac}}\phi_{\tilde{\omega}}, \sqrt{\frac{b-a}{b^2-ac}}\phi_{\tilde{\omega}}\right)$,  and tacking into account the fact that $\phi_{\tilde{\omega}}(x)$ is an  even function on the line, we get $d_{\mathbb{R}}=2d_{\mathbb{R}^+}$.
\end{remark}
\begin{proof} [Proof of Theorem \ref{Theo01}]
 $(i)$ The statement  follows immediately from Lemma \ref{lemma10}.\\
   $(ii)$ We argue by a contradiction. Suppose that \eqref{EeQq} is false. Then, there exist  $\epsilon>0$ and a subsequence $\left\{(u_{n_{k}},v_{n_{k}})\right\}$ such that
\begin{equation*}
\inf_{(\Phi_1,\Phi_2)\in \mathcal{G}(\alpha,\beta)}\left\|(u_{n_{k}},v_{n_{k}})-(\Phi_1,\Phi_2)\right\|_{X}\geq\epsilon\,\,\,\mbox{for any}\,\,\,k\in\mathbb{N}. 
\end{equation*}
But $\left\{(u_{n_{k}},v_{n_{k}})\right\}$  is also a minimizing sequence for $J(\alpha,\beta)$,  hence there exists $(\tilde\Phi_1, \tilde\Phi_2)\in \mathcal{G}(\alpha,\beta)$ such that
\begin{equation*}
\liminf_{k\rightarrow\infty}\left\|(u_{n_{k}},v_{n_{k}})-(\tilde\Phi_1, \tilde\Phi_2)\right\|_{X}=0,
\end{equation*}
which gives a contradiction.
\end{proof}

\subsection{Orbital stability of ground states}\label{subsec_2.3}

This subsection is devoted to the proof of Theorem \ref{theo02} and Corollary \ref{kokoro1}.
First, as an easy consequence of Theorem \ref{Theo01}, we conclude that  the set of minimizers $\mathcal{G}(\alpha,\beta)$ is stable under  the flow generated by system \eqref{Eq}.
\begin{corollary}\label{koro01}
Let $\epsilon>0$. Then there exists $\eta>0$ with the following property:  if $(u_{0},v_{0})\in X$ satisfies 
\begin{equation*}
\inf_{\left(\Phi_1,\Phi_2\right)\in\mathcal{G}(\alpha,\beta)}\left\|(u_{0},v_{0})-\left(\Phi_1,\Phi_2\right)\right\|_{X}<\eta,
\end{equation*}
then the solution $\left(u(t,x),v(t,x)\right)$ of \eqref{Eq} with $\left(u(0,x),v(0,x)\right)=(u_{0},v_{0})$ satisfies
\begin{equation*}
\inf_{\left(\Phi_1,\Phi_2\right)\in\mathcal{G}(\alpha,\beta)}\left\|(u(t,\cdot),v(t,\cdot))-\left(\Phi_1,\Phi_2\right)\right\|_{X}<\epsilon\,\,\,\,\mbox{for all}\,\,\,t\geq0.
\end{equation*}
\end{corollary}
\begin{proof}
By contradiction, we assume that the assertion  is false. Then there exist $\epsilon>0$ and two sequences $\left\{\big(u_{n}(0,x),v_{n}(0,x)\big)\right\}\subset X$ and  $\left\{t_{n}\right\}\subset\mathbb{R}$ such that
\begin{equation}\label{1010}
\inf_{\left(\Phi_1,\Phi_2\right)\in\mathcal{G}(\alpha,\beta)}\left\|\big(u_{n}(0,x),v_{n}(0,x)\big)-\left(\Phi_1,\Phi_2\right)\right\|_{X}<\frac{1}{n}
\end{equation}
and
\begin{equation}\label{1011}
\inf_{\left(\Phi_1,\Phi_2\right)\in\mathcal{G}(\alpha,\beta)}\left\|\big(u_{n}(t_{n},\cdot),v_{n}(t_{n},\cdot)\big)-\left(\Phi_1,\Phi_2\right)\right\|_{X}\geq\epsilon\,\,\,\mbox{for any}\,\,\,n\in\mathbb{N},
\end{equation}
where $\left(u_{n}(t,x),v_{n}(t,x)\right)$ is the solution of \eqref{Eq} with the initial data $\left(u_{n}(0,x),v_{n}(0,x)\right)$. By \eqref{1010} we have that $\left(u_{n}(0,x),v_{n}(0,x)\right)$ converges to an element $\left(\Psi,\Theta\right)\in\mathcal{G}(\alpha,\beta)$ in $X$-norm. Since $Q(\Psi)=\alpha$, $Q(\Theta)=\beta$, and $E(\Psi,\Theta)=J(\alpha,\beta)$, by the conservation laws, we have
\begin{equation*}
Q(u_{n}(t_{n},\cdot))=Q(u_{n}(0,x))\rightarrow\alpha,\,\,\,\,Q(u_{n}(t_{n},\cdot))=Q(u_{n}(0,x))\rightarrow\beta,
\end{equation*}
\begin{equation*}
E\left(u_{n}(t_{n},\cdot),v_{n}(t_{n},\cdot)\right)=E(u_{n}(0,x),v_{n}(0,x))\longrightarrow J(\alpha,\beta)
\end{equation*}
as $n\rightarrow\infty$. Let $\left\{a_{n}\right\}$ and $\left\{b_{n}\right\}$ be such that
\begin{equation*}
Q(a_{n}u_{n}(0,x))=a_n^2Q(u_{n}(0,x))=\alpha\,\,\,\,\mbox{and}\,\,\,\,Q(b_{n}v_{n}(0,x))=b_n^2Q(v_{n}(0,x))=\beta,\qquad n\in\mathbb{N}.
\end{equation*}
Set $\widetilde{u}_{n}(x)=a_{n}u_{n}(x,t_{n})$ and $\widetilde{v}_{n}(x)=b_{n}v_{n}(x,t_{n})$. It is clear that $Q(\widetilde{u}_{n})=\alpha$ and $Q(\widetilde{v}_{n})=\beta$, and since $a_{n},b_{n}\rightarrow1$, we have $\lim\limits_{n\to \infty}E\left(\widetilde{u}_{n},\widetilde{v}_{n}\right)=J(\alpha,\beta)$. Hence $\left\{\left(\widetilde{u}_{n},\widetilde{v}_{n}\right)\right\}$ is a minimizing sequence for $J(\alpha,\beta)$. Thus, by  Theorem \ref{Theo01}-$(ii)$, for $n$ large enough, there exists $\left(\Psi_{n},\Theta_{n}\right)\in\mathcal{G}(\alpha,\beta)$ such that $\left\|\left(\widetilde{u}_{n},\widetilde{v}_{n}\right)-\left(\Psi_{n},\Theta_{n}\right)\right\|_{X}<\frac{\epsilon}{2}$. Since
\begin{equation*}
\begin{split}
&\left\|\left(u_{n}(t_{n},\cdot),v_{n}(t_{n},\cdot)\right)-\left(\Psi_{n},\Theta_{n}\right)\right\|_{X}\\&\leq\left\|\left(u_{n}(t_{n},\cdot),v_{n}(t_{n},\cdot)\right)-\left(\widetilde{u}_{n},\widetilde{v}_{n}\right)\right\|_{X}+\left\|\left(\widetilde{u}_{n},\widetilde{v}_{n}\right)-\left(\Psi_{n},\Theta_{n}\right)\right\|_{X},
\end{split}
\end{equation*}
then, by \eqref{1011}, and $\left\|\left(u_{n}(t_{n},\cdot),v_{n}(t_{n},\cdot)\right)-\left(\widetilde{u}_{n},\widetilde{v}_{n}\right)\right\|_{X}\underset{n\to\infty}{\rightarrow 0}$, we get 

$\epsilon\leq\lim\limits_{n\rightarrow\infty}\left\|\left(u_{n}(t_{n},\cdot),v_{n}(t_{n},\cdot)\right)-\left(\Psi_{n},\Theta_{n}\right)\right\|_{X}\leq\frac{\epsilon}{2}$, which is a contradiction.
\end{proof}

\begin{proof} [Proof of Theorem \ref{theo02}]
Firstly,  by Proposition \ref{lemma4}, for any fixed $\omega>\frac{\gamma^{2}}{N^{2}}$  and $\theta_1,\theta_2\in\mathbb{R}$
\begin{equation*}
\left(e^{i\theta_1}\sqrt{\frac{b-c}{b^{2}-ac}}\Phi_{\omega,\gamma}(\cdot),e^{i\theta_2}\sqrt{\frac{b-a}{b^{2}-ac}}\Phi_{\omega,\gamma}(\cdot)\right)\in\mathcal{G}(\alpha,\beta).
\end{equation*}
 Using Proposition \ref{lemma4} again (see also \eqref{EEQQ1}), for any $(\Phi_1, \Phi_2)\in \mathcal{G}(\alpha, \beta)$ we have
\begin{equation*}
J(\alpha,\beta)=E(\Phi_1,\Phi_2)\geq E_{1}(\Phi_1)+E_{2}(\Phi_2)\geq J_{1}(\alpha)+J_{2}(\beta)=J(\alpha,\beta).
\end{equation*}
Therefore,
\begin{equation*}
\frac{a}{2}\left\|\Phi_1\right\|^{4}_{4}+\frac{c}{2}\left\|\Phi_2\right\|^{4}_{4}+b\int_{\Gamma}\left|\Phi_1\right|^{2}\left|\Phi_2\right|^{2}dx=\frac{b^{2}-ac}{2(b-c)}\left\|\Phi_1\right\|^{4}_{4}+\frac{b^{2}-ac}{2(b-a)}\left\|\Phi_2\right\|^{4}_{4},
\end{equation*}
which implies that
\begin{equation*}
\int_{\Gamma}\left(\sqrt{\frac{b-a}{b-c}}\left|\Phi_1\right|^{2}-\sqrt{\frac{b-c}{b-a}}\left|\Phi_2\right|^{2}\right)^{2}dx=0.
\end{equation*}
Thus,
\begin{equation}\label{20210}
\left(\frac{b-a}{b-c}\right)^{1/4}\left|\Phi_1(x)\right|=\left(\frac{b-c}{b-a}\right)^{1/4}\left|\Phi_2(x)\right|.
\end{equation}
Secondly, notice that for  $\left(\Phi_1,\Phi_2\right)\in\mathcal{G}(\alpha,\beta)$  there exist  Lagrange multipliers  $\omega_{1},\omega_{2}\in\mathbb{R}$ such that $\left(\Phi_1,\Phi_2\right)$ is a solution of 
 \begin{equation*}
\begin{cases} 
-\Delta_{\gamma}\Phi_1+\omega_{1}\Phi_1=a\left|\Phi_1\right|^{2}\Phi_1+b\left|\Phi_2\right|^{2}\Phi_1\\
-\Delta_{\gamma}\Phi_2+\omega_{2}\Phi_2=b\left|\Phi_1\right|^{2}\Phi_2+c\left|\Phi_2\right|^{2}\Phi_2
\end{cases} 
\end{equation*}
in a weak sense. Moreover, using \eqref{20210}, we  simplify the system 
\begin{equation}\label{3010}
\left\{
\begin{array}{c}
-\Delta_{\gamma}\Phi_1+\omega_{1}\Phi_1=\frac{b^{2}-ac}{b-c}\left|\Phi_1\right|^{2}\Phi_1\\
-\Delta_{\gamma}\Phi_2+\omega_{2}\Phi_2=\frac{b^{2}-ac}{b-a}\left|\Phi_2\right|^{2}\Phi_2.
\end{array}\right.
\end{equation}
Following the arguments from  the proofs of  \cite[Lemmas 25 and 26]{RJJ}, we can show that the only pair of $L^{2}$-solution to \eqref{3010} is given by
\begin{equation*}
\Phi_1(x)=e^{i\theta_1}\sqrt{\frac{2\omega_{1}(b-c)}{b^{2}-ac}}\,\Phi_{\omega_{1},\gamma}
\end{equation*}
and
\begin{equation*}
\Phi_2(x)=e^{i\theta_2}\sqrt{\frac{2\omega_{2}(b-a)}{b^{2}-ac}}\,\Phi_{\omega_{2},\gamma},
\end{equation*}
where $\theta_1,\theta_2\in\mathbb{R}$. From \eqref{20210} we get 
$$\left|\Phi_{\omega_{2},\gamma}(x)\right|=\left|\sqrt{\frac{\omega_1}{\omega_2}}\Phi_{\omega_{1},\gamma}(x)\right|.$$
Finally, since
\begin{equation*}
\left\|\Phi_1\right\|^{2}_{2}=\alpha=2\frac{b-c}{b^{2}-ac}(N\sqrt{\omega}-\gamma)\,\,\,\mbox{and}\,\,\, \left\|\Phi_2\right\|^{2}_{2}=\beta=2\frac{b-a}{b^{2}-ac}(N\sqrt{\omega}-\gamma), 
\end{equation*}
and using \eqref{20210}, it follows that $\omega_{1}=\omega_{2}=\omega$. This completes the proof.
\end{proof}

To end this section notice that  Corollary \ref{kokoro1} follows from  Corollary \ref{koro01} and Theorem \ref{theo02}.

%%%%%%%%%%%%%%%%%%%%%%%%%%%%%%%%%%%%%%%%%%%%%%%%%%%%%%%%%%%%%%%%%%%%%%%%%%%%%%%%%%%%%%%%%%%%%%%%%%%%%%%%%%%%%%%%5
\section{Orbital stability of standing waves: case of generalized nonlinearity}\label{sec_3}
In this section we study orbital stability of standing waves of the general system 
\begin{equation}\label{syst_geral}
\left\{\begin{array}{l}
i \partial_{t} u(t, x)+\Delta_{\gamma} u(t, x)+\left(a|u(t, x)|^{q-2}+b|v(t, x)|^{p}|u(t,x)|^{p-2}\right) u(t, x)=0 \\
i \partial_{t} v(t, x)+\Delta_{\gamma} v(t, x)+\left(c|v(t, x)|^{r-2}+b|u(t,x)|^p|v(t, x)|^{p-2}\right) v(t, x)=0\\
(u(0, x), v(0, x))=\left(u_{0}(x), v_{0}(x)\right).
\end{array}\right.
\end{equation}
We assume that $a,b, c\in\mathbb{R}$ and $2<q,r, 2p$.   The  system is  locally well posed (the proof follows analogously to \cite[Theorem 4.10.1]{Caz03}, see also \cite[Section 2]{AQFF}). 
In particular,  the conserved  energy  is given by
\begin{equation}\label{energy_geral}
E(u, v)=\|u'\|_2^2+\|v'\|_2^2-\gamma\left(|u_{1}(0)|^2+|v_{1}(0)|^2\right)-G(u,v),
\end{equation}
where 
$$G(u,v)= \frac{2a}{q}\|u\|^{q}_q+\frac{2c}{r}\|v\|^{r}_r+\frac{2b}{p}\|uv\|^{p}_p.$$ 
\begin{remark}\label{general_system}
Notice that for $q,r, 2p<6$ the global well-posedness holds. 
Due to the conservation of masses \eqref{masses} and energy, we have 
\begin{equation}\label{X_norm}
\begin{split}
\|(u,v)\|_X^2&=E(u,v)+\gamma\left(|u_1(0)|^2+|v_1(0)|^2\right)+G(u,v)+\|u\|_2^2+\|v\|_2^2\\&=\|(u_0,v_0)\|_X^2-\gamma\left(|u_{01}(0)|^2+|v_{01}(0)|^2\right)-G(u_0,v_0)+G(u,v).
\end{split}
\end{equation}
Observe that 
\begin{equation}\label{G}
\begin{split}
G(u,v)&\leq C_1\|u'\|_2^{\frac{q-2}{2}}\|u\|_2^{\frac{q+2}{2}}+C_2\|v'\|_2^{\frac{r-2}{2}}\|v\|_2^{\frac{r+2}{2}} +\frac{2b}{p}\|u\|_{2p}^p\|v\|_{2p}^p \\&\leq  C_1\|u'\|_2^{\frac{q-2}{2}}\|u\|_2^{\frac{q+2}{2}}+C_2\|v'\|_2^{\frac{r-2}{2}}\|v\|_2^{\frac{r+2}{2}} +\frac{b}{p}\left(\|u\|_{2p}^{2p}+\|v\|_{2p}^{2p}\right)\\&\leq \varepsilon\left(\|u'\|_2^2+\|v'\|_2^2\right) +C_3(\varepsilon)\left(\|u\|_2^{\frac{2(q+2)}{6-q}}+\|v\|_2^{\frac{2(r+2)}{6-r}}+\|u\|_2^{\frac{2(2p+2)}{6-2p}}+\|v\|_2^{\frac{2(2p+2)}{6-2p}}\right). 
\end{split}    
\end{equation}
 The above estimate is induced by  Gagliardo-Nirenberg inequality \eqref{GNIq}  and the Young inequality $fg \leq \varepsilon f^{l}+C_{\varepsilon} g^{l^{\prime}}, \frac{1}{l}+\frac{1}{l^{\prime}}=$ $1, l, l^{\prime}>1,\,  f, g \geq 0$. Observe that the key point is that $l=\frac{4}{q-2}>1$ for $2<q<6$ (analogously for $r$ and $2p$).
From \eqref{Gagli1}, \eqref{X_norm}, and \eqref{G} we get 
$$\|(u,v)\|_X^2\leq \frac{1}{1-\varepsilon}\Big\{\|(u_0,v_0)\|_X^2-G(u_0,v_0)+C(\|u_0\|_2^2,\|v_0\|_2^2)\Big\},$$
i.e.  the norm  $\|(u,v)\|_X^2$ is controlled by the constant independent on time. 
Finally, repeating the proof of \cite[Theorem 3.4.1]{Caz03} (starting from formula (3.4.2)), we obtain the  global existence.
\end{remark}

\subsection{Stability of  one component standing waves }\label{subsec_3.1} 
We consider  the standing waves  of the simplest  form $(e^{i\omega t}\Phi_1,0)$ and $(0, e^{i\omega t}\Phi_2)$. Throughout this section, we assume additionally $p\geq 2$.
Observe that  $(\Phi_1,0), (0,\Phi_2)$ are critical points of the  functional 
\begin{equation*}\label{action}
S_\omega(u,v)=\frac{1}{2}\left\{E(u,v)+\omega( \|u\|_2^2+\|v\|_2^2)\right\},\end{equation*}
that is, they satisfy equation \eqref{phi_12}.
The description of the solutions to the first equation \eqref{phi_12} was given in \cite[Theorem 4]{AQFF} (to obtain description to the second equation we need to replace $q$ by $r$, and $a$ by $c$):
\begin{theorem}\label{description} Let $[s]$ denote the integer part of $s \in \mathbb{R}, \gamma \neq 0$, and $a>0$. Then first equation in \eqref{phi_12} has $\left[\frac{N-1}{2}\right]+1$ (up to permutations of the edges of $\Gamma$ and rotation) vector solutions $\Phi_{k}^{\gamma}=\left(\varphi_{k, j}^{\gamma}\right)_{j=1}^{N}, k=0, \ldots,\left[\frac{N-1}{2}\right]$, which are given by
\begin{equation*}
\begin{aligned}
\varphi_{k, j}^{\gamma}(x)=\left\{\begin{array}{ll}
{\left[\frac{q \omega}{2a} \operatorname{sech}^{2}\left(\frac{(q-2) \sqrt{\omega}}{2} x-a_{k}\right)\right]^{\frac{1}{q-2}},} & j=1, \ldots, k ; \\
{\left[\frac{q \omega}{2a} \operatorname{sech}^{2}\left(\frac{(q-2) \sqrt{\omega}}{2} x+a_{k}\right)\right]^{\frac{1}{q-2}},} & j=k+1, \ldots, N,
\end{array}\right. \\
\text { where } a_{k}=\tanh ^{-1}\left(\frac{\gamma}{(N-2k) \sqrt{\omega}}\right), \text { and } \omega>\frac{\gamma^{2}}{(N-2 k)^{2}}.
\end{aligned}
\end{equation*}
\end{theorem}
Below we deal with two types of instability: orbital and spectral. General definition of the orbital stability for a Hamiltonian system invariant under the action of some Lie group can be found in \cite[Section 2]{GrilSha90}. 
The definition of the spectral instability involves a  linearization of \eqref{syst_geral}  around the profile of the
standing wave. After making necessary technical steps we give precise Definition \ref{spec_instub}.

First,  we linearize system \eqref{syst_geral} around $(\Phi_k^\gamma,0)$. We observe that  system \eqref{syst_geral} can be written in the form 
\begin{equation*}\label{hamilt}\frac{d}{d t}\left(\begin{array}{l}
u \\
v
\end{array}\right)=\frac{1}{2}\mathcal{J} E^{\prime}[u, v],\quad \mathcal{J}=\left(\begin{array}{ll}
-i & 0\\
0 & -i
\end{array}\right),\end{equation*} and put 
$$(u(t), v(t))=e^{i \omega t}\left\{(\Phi_k^\gamma,0)+(h_1(t),h_2(t))\right\}.$$ Since  $S_{\omega}^{\prime}\left(\Phi_k^\gamma,0\right)=0$,
we get 
$$\frac{d}{d t}\left(\begin{array}{l}
h_1\\
h_2
\end{array}\right)=\mathcal{J} S_{\omega}^{\prime \prime}\left(\Phi_k^\gamma,0\right)\left(\begin{array}{l}
h_1 \\
h_2
\end{array}\right) +O\left(\|(h_1,h_2)\|_{X}^{2}\right), $$
where
\begin{align*}
 &S_{\omega}^{\prime \prime}\left(\Phi_k^\gamma,0\right)=  \left(\begin{array}{ll}
S_1 & 0\\
0 & S_2
\end{array}\right), \quad S_1 h_1=-\widetilde{\Delta}_{\gamma} h_1+\omega h_1-a\left(\Phi_k^\gamma\right)^{q-2} h_1-a(q-2)\left(\Phi_k^\gamma\right)^{q-2} \operatorname{Re}(h_1),\\
&S_2h_2=\left\{\begin{array}{c}
     -\widetilde{\Delta}_{\gamma} h_2+\omega h_2-b\left(\Phi_k^\gamma\right)^2 h_2, \,\, p=2 \\
      -\widetilde{\Delta}_{\gamma} h_2+\omega h_2,\,\, p>2.\qquad\qquad\quad
\end{array}\right.
\end{align*}
Here $-\widetilde{\Delta}_{\gamma}: H^{1}(\Gamma) \rightarrow\left(H^{1}(\Gamma)\right)^*$ is the  unique bounded operator associated with the bounded on $H^1(\Gamma)$ bilinear form $t_{\gamma}\left(u_{1}, u_{2}\right)=\left(u_{1}^{\prime}, u_{2}^{\prime}\right)_{2}-\gamma \operatorname{Re}\left(u_{11}(0) \overline{u_{21}(0)}\right)$.
By the Representation Theorem \cite[Chapter VI, Theorem 2.1]{Kat66}, we can associate with the bilinear form 
\begin{equation*}
\begin{split}
b_{\gamma}\big((h_1, h_2), (z_1,z_2)\big)=\left\langle S_{\omega}^{\prime \prime}\left(\Phi_k^\gamma,0\right) (h_1,h_2),(z_1,z_2)\right\rangle_{X^*\times X}
\end{split}
\end{equation*}
 self-adjoint in $L^{2}(\Gamma)\times L^{2}(\Gamma)$ (with the real scalar product) operator
\begin{align*}
 &L=  \left(\begin{array}{ll}
L_1 & 0\\
0 & L_2
\end{array}\right), \quad L_1 h_1=- h''_1+\omega h_1-a\left(\Phi_k^\gamma\right)^{q-2} h_1-a(q-2)\left(\Phi_k^\gamma\right)^{q-2} \operatorname{Re}(h_1),\\
&L_2h_2=\left\{\begin{array}{c}
     -h''_2+\omega h_2-b\left(\Phi_k^\gamma\right)^2 h_2, \,\, p=2 \\
      -h''_2+\omega h_2,\,\, p>2.\qquad\qquad\quad
\end{array}\right.,\quad \dom(L_1)=\dom(L_2)=\dom(\Delta_\gamma).
\end{align*}
Putting $h_j=h_j^R+ih_j^I, j=1,2,$ and identifying $L_{\mathbb{C}}^{2}(\Gamma)$ with $L_{\mathbb{R}}^{2}(\Gamma) \oplus L_{\mathbb{R}}^{2}(\Gamma)$, we conclude
\begin{equation}\label{equivalence}\mathcal J L\quad \Longleftrightarrow \left(\begin{array}{cccc}
   0  & I_{L^2_{\mathbb{R}}} & 0&0\\
   -I_{L^2_{\mathbb{R}}}  & 0 & 0&0\\
   0 & 0 &0& I_{L^2_{\mathbb{R}}}\\
   0 &0 & -I_{L^2_{\mathbb{R}}}&0
\end{array}\right)\left(\begin{array}{cccc}
   L_1^R & 0 & 0&0\\
   0  & L_1^I & 0&0\\
   0 & 0 & L_2^R & 0\\
   0 &0 & 0& L_2^I
\end{array}\right),\end{equation} that is, 
$$\mathcal J L\quad \Longleftrightarrow \left(\begin{array}{cc}
    J & 0  \\
    0 & J
\end{array}\right)\left(\begin{array}{cc}
    L_1 & 0  \\
    0 & L_2
\end{array}\right),$$
where $J=\left(\begin{array}{cc}
    0 & I_{L^2_{\mathbb{R}}}  \\
    -I_{L^2_{\mathbb{R}}} & 0
\end{array}\right)$ and 
\begin{align*}
&L_1^R=-\frac{d^2}{dx^2}+\omega-a(q-1)(\Phi_k^\gamma)^{q-2}, \quad L_1^I=-\frac{d^2}{dx^2}+\omega-a(\Phi_k^\gamma)^{q-2},\\
&L_2^R=L_2^I=L_2, \quad  \dom(L_j^R)=\dom(L_j^I)=\dom(\Delta_\gamma),\,\, j=1,2.
\end{align*}
\begin{definition}\label{spec_instub} The standing wave $e^{i\omega t }(\Phi^\gamma_k, 0)$ is said to be spectrally unstable if there exist $\lambda$ with $\operatorname{Re} \lambda>0$ and $\vec{h}=(h_1, h_2)\in \operatorname{dom}\left(\Delta_{\gamma}\right)\times \dom(\Delta_\gamma)$ such that
$$\mathcal{J}L\vec{h}=\lambda \vec{h}.$$
\end{definition}
The notion of spectral instability is particularly important since frequently its presence leads to nonlinear instability.
\begin{remark}
$(i)$ In a view  of  \eqref{equivalence} the above definition can be rewritten as $$
\left(\begin{array}{cc}
    J & 0  \\
    0 & J
\end{array}\right)\left(\begin{array}{cc}
    L_1 & 0  \\
    0 & L_2
\end{array}\right)\left(\begin{array}{c}
    h_1 \\
    h_2
\end{array}\right)=\lambda\left(\begin{array}{c}
    h_1 \\
    h_2
\end{array}\right) .
$$

$(ii)$ Notice that spectral instability implies that $(0,0)$ is unstable solution to the linearized equation
$$
\frac{d}{d t}\left(\begin{array}{l}
u(t) \\
v(t)
\end{array}\right)=\mathcal{J} S_{\omega}^{\prime \prime}\left(\Phi_{k}^{\gamma}, 0\right)\left(\begin{array}{l}
u(t) \\
v(t)
\end{array}\right)
$$
in the sense of Lyapunov.
\end{remark}
Below we list spectral properties of the operators $L_1^R, L_1^I$ proved in  \cite[Proposition 6.1]{AQFF} and \cite[Theorem 3.1]{Kai19}:  let $k=0, \ldots,\left[\frac{N-1}{2}\right]$, then 
\begin{itemize}
    \item[$a)$] $\ker(L_1^I)=\Span\{\Phi_k^\gamma\}$ and     $\ker(L_1^I)=\{0\}$;
    \item[$b)$] $L_1^I\geq 0$ and  $n\left(L_1^R\right)=\left\{\begin{array}{l}k+1 \quad \text { if } \quad \gamma>0, \\ N-k \quad \text { if } \quad \gamma<0\end{array}\right.;$
    \item[$c)$] $\sigma_{\ess}(L_1^R)=\sigma_{\ess}(L_1^I)=[\omega, \infty).$
\end{itemize}
We use these properties to show the principal result of this subsection.
\begin{theorem}\label{spec_inst_1_hump}
Let  $a>0$, $q>2$, and $\Phi^\gamma_k$ be defined in Theorem \ref{description}.

$(i)$\, Assume  $p\geq 2$.  If either $\gamma>0, k\geq 1$ or $\gamma<0, k\geq 0$, then   $e^{i\omega t }(\Phi^\gamma_k, 0)$ is spectrally unstable. If additionally  $p,q>3$, then we have orbital instability.

$(ii)$\, Assume $p>2$. If $\gamma>0, k=0$, then $e^{i\omega t }(\Phi^\gamma_0, 0)$ is orbitally stable for $2<q\leq 6$. Moreover, for $q>6$ there exists $\omega_1>\frac{\gamma^2}{N^2}$ such that     $e^{i\omega t }(\Phi^\gamma_0, 0)$ is orbitally stable as $\omega \in (\frac{\gamma^2}{N^2}, \omega_1)$, while it is orbitally unstable as $\omega>\omega_1$.
\end{theorem}
\begin{proof}
$(i)$\, From \cite[Theorem 1.2]{Grill88} one concludes that the operator $JL_1$ has positive eigenvalue for $\gamma>0, k\geq 1,$ and $\gamma<0, k\geq 0$. 
Indeed, $JL_1\vec{w}=\lambda\vec{w}, \vec{w}=(w_1,w_2)$, is equivalent to  
\begin{equation}\label{lambda_1}
\left(\begin{array}{cc}
0 &  L^I_1 \\
- L^R_1 & 0
\end{array}\right)\left(\begin{array}{l}
w_{1} \\
w_{2}
\end{array}\right)=\lambda\left(\begin{array}{c}
w_{1} \\
w_{2}
\end{array}\right),
\end{equation}
or  $\left\{\begin{array}{c}
 L^I_1 w_{2}=\lambda w_{1} \\
 L^R_1 w_{1}=-\lambda w_{2}
\end{array}\right.$. Since $\ran(L^I_1) \perp \ker(L^I_1)$, we get $w_{1} \in \ker(L^I_1)^{\perp}$. Hence $w_{2}=\lambda(L^I_1)^{-1} w_{1}$ and $L^R_1 w_{1}+\lambda^{2}(L^I_1)^{-1} w_{1}=0 .$  

Consider the problem 
$$P L^R_1 w_{1}+\lambda^{2}P(L^I_1)^{-1} w_{1}=0,$$
 where $P$  is the orthogonal projection onto $$\ker( L^I_1)^{\perp}=\left\{(v_1,v_2)\in L^{2}(\Gamma)\times L^{2}(\Gamma) :\big((v_1,v_2), (\Phi^\gamma_k,0)\big)_{2}=0\right\}$$ (here we assume that $L^{2}(\Gamma)$ is endowed with the usual complex inner product).
The projection $P$ serves to fit the problem into the Hilbert space $\ker(L^I_1)^{\perp} .$
\cite[Theorem 1.2]{Grill88}  states that the number $I(L^R_1,  L^I_1)$ of positive  $\lambda$ satisfying \eqref{lambda_1} is estimated by
\begin{equation}\label{estim_lambda_1}
n(P L^R_1)-n(P(L^I_1)^{-1}) \leq I(L^R_1,  L^I_1).
\end{equation}
Using spectral properties $a),b)$ mentioned before the theorem,  formula \eqref{estim_lambda_1} yields 
$$I(L^R_1,  L^I_1)\geq\left\{\begin{array}{l}k,  \quad \gamma>0, \\
N-k-1,\quad \gamma<0.\end{array}\right.$$

By \eqref{equivalence}, there exists $\lambda>0$ satisfying Definition \ref{spec_instub}  under the assumptions in $(i)$, that is, $(\Phi^\gamma_0, 0)$ is spectrally unstable. Using the fact that the nonlinearity $$F_{p,q,r}(u,v)=(a|u|^{q-2}u+b|v|^p|u|^{p-2}u, c|v|^{r-2}v+b|u|^p|v|^{p-2}v)$$ is of $C^2$ class for $p, q>3$, we conclude  that mapping data-solution associated to \eqref{syst_geral} is of class $C^2$ (for instance, see \cite[Step 4 of proof of Proposition 2.1]{Gol21}).  Finally,  to imply the orbital instability from the spectral one, the approach by \cite{HenPer82} (see Theorem 2 and \cite[Theorem 5.8, Corollary 5.9]{Gol21}) can be used.

$(ii)$\, By $L_2\geq \omega$, we deduce  $n(L)=n(L_1)=1$. Observe that $S''_{\omega}(\Phi^\gamma_0, 0)$  satisfies \cite[Condition (G)]{Stu08}. Indeed, for $\vec{v}=(v_1, v_2)\in X$ 
\begin{equation*}
\begin{split}
 &\langle S''_{\omega}(\Phi^\gamma_0, 0)\vec{v}, \vec{v}\rangle_{X^*\times X}=\int_{\Gamma}\Big\{|v'_1|^2+|v'_2|^2+\omega\Big(|v_1|^2+|v_2|^2\Big)\Big\}dx
 \\& -\int_{\Gamma}\Big\{a(\Phi_0^\gamma)^{q-2}|v_1|^2+a(q-2)(\Phi_0^\gamma)^{q-2}(\re v_1)^2\Big\}dx
 -\gamma\Big(|v_{11}(0)|^2+|v_{21}(0)|^2\Big)\\&
 \geq \min\{\frac{1}{2}, \omega\}\|\vec{v}\|^2_X-\Big((q-1)aM^{q-2}+\frac{2\gamma^2}{N^2}\Big)\|\vec{v}\|^2_2,
 \end{split}    
\end{equation*}
where $M=\|\Phi_0^\gamma\|_\infty$ and we have used
$$-\gamma\Big(|v_{11}(0)|^2+|v_{21}(0)|^2\Big)+\frac{1}{2}\int_{\Gamma}\Big\{|v'_1|^2+|v'_2|^2\Big\}dx\geq -\frac{2\gamma^2}{N^2}\|\vec{v}\|_2^2.$$
Let $\mathcal{R}: X\rightarrow X^*$ be  the Riesz isomorphism. Thus, from \cite[Lemma 5.4]{Stu08} we conclude $\sigma_{\ess}(\mathcal{R}^{-1}S''_{\omega}(\Phi^\gamma_0,0))=\sigma_{\ess}(L)=[\omega, \infty)$, $\ker(\mathcal{R}^{-1}S''_{\omega}(\Phi^\gamma_0,0))=\ker(L)=\Span\{i(\Phi^\gamma_0,0)\}$, and $n(\mathcal{R}^{-1}S''_{\omega}(\Phi^\gamma_0,0))=n(L)=1$.  

  Denote $d''(\omega)=\partial^2_\omega(S_\omega(\Phi^\gamma_0,0))$. By \cite[Proposition 3.19-$(i)$]{AngGol18}, we get  $d''(\omega)>0$ for $2<q\leq 6$, and for $q>6$ there exists $\omega_1$ such that $d''(\omega)>0$ as $\omega\in(\frac{\gamma^2}{N^2}, \omega_1)$, and $d''(\omega)<0$ as $\omega>\omega_1$. Finally, the result follows applying \cite[Theorem 3]{GrilSha87}.
\end{proof}
\begin{remark}
$(i)$\,
Observe that for $p=2, q=4, \gamma>0, k=0$, and $b<a$ we have $L_2>L_1^I$, then $\sigma(L_2)=[\omega, \infty).$ Therefore, $n(\mathcal{R}^{-1}S''_\omega(\Phi_0^\gamma,0))=1$. As in the proof of item $(ii)$, we have $d''(\omega)>0$, then the orbital stability of $e^{i\omega t}(\Phi^\gamma_0, 0)$ holds.

$(ii)$\, All the above results on the stability/instability of $(\Phi^\gamma_k, 0)$ can be repeated for the profile $(0, \Phi_2)=(0, \Phi^\gamma_k)$ with $q$ replaced by $r$ and $a$ replaced by $c$.
\end{remark}

\subsection{Stability of  two component standing waves}\label{subsec_3.2}
In this subsection we study stability properties of the standing waves of  system \eqref{syst_geral} of the form $(e^{i\omega t}\Phi(x), e^{i\omega t}\Phi(x))$  in the case $q=r=2p-2, a=c=1$:
\begin{equation}\label{syst_particular}
\left\{\begin{array}{l}
i \partial_{t} u(t, x)+\Delta_{\gamma} u(t, x)+\left(|u(t, x)|^{2p-2}+b|v(t, x)|^{p}|u(t,x)|^{p-2}\right) u(t, x)=0 \\
i \partial_{t} v(t, x)+\Delta_{\gamma} v(t, x)+\left(|v(t, x)|^{2p-2}+b|u(t,x)|^p|v(t, x)|^{p-2}\right) v(t, x)=0.
\end{array}\right.
\end{equation}
 It is easily seen that  $\Phi(x)$ solves
\begin{equation}\label{profile_multi}
-\Delta_\gamma \Phi+\omega \Phi-(b+1) \Phi^{2 p-1}=0.
\end{equation}
The description of solutions to \eqref{profile_multi} is given by Theorem \ref{description} (one just needs to substitute $q$ by $2p$ and $a$ by $b+1$).
Namely, for $b>-1$ we have $\left[\frac{N-1}{2}\right]+1$ solutions $\Phi_{k}^{\gamma}=\left(\varphi_{k, j}^{\gamma}\right)_{j=1}^{N}$, $k=0, \ldots,\left[\frac{N-1}{2}\right],$ of the form 
\begin{equation}\label{hump}
\begin{aligned}
\varphi_{k, j}^{\gamma}(x)=\left\{\begin{array}{ll}
{\left[\frac{p \omega}{(b+1)} \operatorname{sech}^{2}\Big((p-1) \sqrt{\omega} x-a_{k}\Big)\right]^{\frac{1}{2p-2}},} & j=1, \ldots, k ; \\
{\left[\frac{p \omega}{(b+1)} \operatorname{sech}^{2}\Big((p-1) \sqrt{\omega} x+a_{k}\Big)\right]^{\frac{1}{2p-2}},} & j=k+1, \ldots, N,
\end{array}\right. \\
\text { where } a_{k}=\tanh ^{-1}\left(\frac{\gamma}{(N-2k) \sqrt{\omega}}\right), \text { and } \omega>\frac{\gamma^{2}}{(N-2 k)^{2}}.
\end{aligned}
\end{equation} 
The main results of this subsection are  the following two theorems. The first one deals with arbitrary $k$ and the second one contains stability results for the particular case of symmetric tail-profile $(\Phi_0^\gamma, \Phi_0^\gamma)$.
\begin{theorem}\label{main_1}
Let $p\geq 2, b>-1$, and  $(\Phi^\gamma_k, \Phi^\gamma_k)$ be defined by \eqref{hump}. If $\gamma>0, k\geq 2$ or $\gamma<0, N-k\geq 3$, then  the standing wave $e^{i\omega t}(\Phi^\gamma_k, \Phi^\gamma_k)$  is spectrally unstable. Moreover, for $p>3$ the orbital instability holds. 

If additionally $0<b<p-1$, then the above assertions hold for $\gamma<0, k\geq 0$ and $\gamma>0, k=1$.
\end{theorem}
\begin{theorem}\label{main_2}
Let $\gamma>0, k=0,$ then for  $2\leq p<3$ and $0<b\neq p-1$ the standing wave $e^{i\omega t}(\Phi^\gamma_0, \Phi^\gamma_0)$ is orbitally stable, while for $p>3$ and $b>p-1$ the standing wave $e^{i\omega t}(\Phi^\gamma_0, \Phi^\gamma_0)$ is orbitally unstable.
\end{theorem}
As in the previous case, we linearize \eqref{syst_particular} at $(\Phi^\gamma_k, \Phi^\gamma_k)$. Notice that the profile $(\Phi_1(x), \Phi_2(x))$ of the standing wave  $(e^{i\omega_1 t}\Phi_1(x), e^{i\omega_2 t}\Phi_2(x))$ is the critical point of the action functional given by 
$$S_{\omega_1, \omega_2}(u,v)=\frac{1}{2}\left\{E(u,v)+\omega_1\|u\|_2^2+\omega_2\|v\|_2^2\right\},$$ where $E$ is defined by \eqref{energy_geral}.
We obtain for $\vec{h}=(h_1,h_2)\in X$
 \begin{equation*}
 \begin{split}
 &S_{\omega,\omega}^{\prime \prime}\left(\Phi_k^\gamma,\Phi_k^\gamma\right)\vec{h}=\left(\begin{array}{l}
\tilde{S_1}\vec{h} \\
\tilde{S_2}\vec{h}
\end{array}\right),\,\,\tilde{S_1}\vec{h}=\tilde{{l}_\gamma}h_1-bp(\Phi^\gamma_k)^{2p-2}\re(h_2), \,\, \tilde{S_2}\vec{h}=\tilde{{l}_\gamma}h_2-bp(\Phi^\gamma_k)^{2p-2}\re(h_1),  \\
&\tilde{l}_\gamma=-\tilde{\Delta}_\gamma +\Big(\omega-(1+b)(\Phi^\gamma_k)^{2p-2}\Big)-\Big(\big(2p-2+b(p-2)\big)(\Phi^\gamma_k)^{2p-2}\Big)\re(\cdot). \end{split}\end{equation*}
Analogously to the previous case,  we associate with the  bilinear form \\
$\left\langle S_{\omega,\omega}^{\prime \prime}\left(\Phi_{k}^{\gamma}, \Phi_{k}^{\gamma}\right)\left(h_{1}, h_{2}\right),\left(z_{1}, z_{2}\right)\right\rangle_{X^{*} \times X}$ the self-adjoint in $L^2(\Gamma)\times L^2(\Gamma)$ operator  $\tilde{L}$. Putting $h_j=h_j^R+ih_j^I, j=1,2,$ we obtain 
\begin{equation}\label{4.13}
\begin{split}&\tilde L=\tilde L^R\vec{h}^R+i \tilde L^I\vec{h}^I, \quad \tilde L^I=\left(\begin{array}{cc}
    \tilde L^I_1 & 0 \\
   0  & \tilde L^I_2
\end{array}\right),\,\, \tilde L^I_{j}=-\frac{d^2}{dx^2}+\omega-(1+b)(\Phi^\gamma_k)^{2p-2}, \\
&\tilde L^R=\Small{\left(\begin{array}{cc}
  -\frac{d^2}{dx^2}+\omega-\big(2p-1+b(p-1)\big)(\Phi^\gamma_k)^{2p-2}   & -bp(\Phi^\gamma_k)^{2p-2}  \\
   -bp(\Phi^\gamma_k)^{2p-2}  &  -\frac{d^2}{dx^2}+\omega-\big(2p-1+b(p-1)\big)(\Phi^\gamma_k)^{2p-2}
\end{array}\right)}, 
\end{split}\end{equation}
where $\dom(\tilde L^I)=\dom(\tilde L^R)=\dom(\Delta_\gamma)\times \dom(\Delta_\gamma).$
It is easily seen that 
\begin{equation}\label{equivalence_hump}
\mathcal{J}\tilde L\,\,  \Longleftrightarrow\,\, \left(\begin{array}{cc}
    0 & I_{L^2_{\mathbb{R}}\times L^2_{\mathbb{R}}} \\
   -I_{L^2_{\mathbb{R}}\times L^2_{\mathbb{R}}}  & 0
\end{array}\right)\left(\begin{array}{cc}
    \tilde L^R & 0 \\
   0  & \tilde L^I
\end{array}\right)=\left(\begin{array}{cc}
   0 &  \tilde L^I \\
    -\tilde L^R & 0
\end{array}\right).        
\end{equation}
\begin{remark}
The definition of the spectral instability is analogous to Definition \ref{spec_instub}. One just needs to substitute $L$ by $\tilde L$. By \eqref{equivalence_hump}, it means the existence of $\lambda$ with $\re \lambda>0$ and $(h_1,h_2)\in \dom(\Delta_\gamma)\times \dom(\Delta_\gamma)$ such that 
$$\left(\begin{array}{cc}
   0 &  \tilde L^I \\
    -\tilde L^R & 0
\end{array}\right)\left(\begin{array}{c}
   h_1\\
   h_2
\end{array}\right)=\lambda \left(\begin{array}{c}
   h_1\\
   h_2
\end{array}\right).$$
\end{remark}
Below we list some properties of the operators $\tilde L^I$ and $\tilde L^R$ (see \cite[Proposition 6.1]{AQFF}):
\begin{itemize}
    \item[$a)$] $\ker(\tilde L^I)=\Span\{(\Phi^\gamma_k, 0), (0, \Phi^\gamma_k)\}$;   
    \item[$b)$] $\tilde L^I\geq 0$; 
    \item[$c)$] $\sigma_{\ess}(\tilde L^I)=\sigma_{\ess}(\tilde L^R)=[\omega, \infty)$.
    \end{itemize}
    \begin{remark}
    The property $\sigma_{\ess}(\tilde L^R)=[\omega, \infty)$ follows from the representation 
    \begin{align*}
        \tilde L^R&=\Small{\left(\begin{array}{cc}
  -\frac{d^2}{dx^2}+\omega-\big(2p-1+b(p-1)\big)(\Phi^\gamma_k)^{2p-2}    & 0  \\
  0  &  -\frac{d^2}{dx^2}+\omega-\big(2p-1+b(p-1)\big)(\Phi^\gamma_k)^{2p-2}
\end{array}\right)}\\&+\Small{\left(\begin{array}{cc}
 0  & -bp(\Phi^\gamma_k)^{2p-2}  \\
  -bp(\Phi^\gamma_k)^{2p-2}  &  0
\end{array}\right)} \end{align*} and the fact that $ -bp(\Phi^\gamma_k)^{2p-2} $ is relatively $-\Delta_\gamma$-compact. 

\end{remark}
Below we count $n(\tilde L^R)$.
\begin{proposition}\label{count_n}
Let $p>-1$, then $n(\tilde L^R)\geq \left\{\begin{array}{c}
     k+1,\,\, \gamma>0 \\
     N-k,\,\, \gamma<0
\end{array}\right.$ with $k=0, \ldots,\left[\frac{N-1}{2}\right]$. Moreover, if additionally

$a)$\, $0<b<p-1$, then $n(\tilde L^R)\geq \left\{\begin{array}{c}
     k+2,\,\, \gamma>0 \\
     N-k+1,\,\, \gamma<0
\end{array}\right.;$

$b)$\, $b>p-1$, then $n(\tilde L^R)= \left\{\begin{array}{c}
     k+1,\,\, \gamma>0 \\
     N-k,\,\, \gamma<0.
\end{array}\right.$
\end{proposition}
To prove the above proposition, we suppose that $\lambda\in \sigma_p(\tilde L^R)$ and $(h_1, h_2)$ is the corresponding eigenvector.  Then denoting $h_+=h_1+h_2$ and $h_-=h_1-h_2, $ we get
\begin{equation}\label{decomposition}
  \begin{split}
    &\left\{\begin{array}{c}
      \tilde L^R_+h_+=\lambda h_+    \\
       \tilde L^R_-h_-=\lambda h_-  
    \end{array}\right.,\quad  \tilde L^R_+=-\frac{d^2}{dx^2}+\omega-(2p-1)(b+1)(\Phi_k^\gamma)^{2p-2},   \\
    & \tilde L^R_-=-\frac{d^2}{dx^2}+\omega-(2p-b-1)(\Phi_k^\gamma)^{2p-2},\quad \dom(\tilde L^R_+)= \dom(\tilde L^R_-)=\dom(\Delta_\gamma). 
  \end{split}  
\end{equation}
From \eqref{decomposition} we conclude 
\begin{equation}\label{decomp_ker}
    (h_1, h_2)\in \ker(\tilde L^R) \quad\Longleftrightarrow\quad (h_+, h_-)\in \ker\left(\begin{array}{cc}
     \tilde L^R_+  & 0   \\
      0   & \tilde L^R_-
    \end{array}\right)
\end{equation}
and
\begin{equation}\label{decomp_n}
  n(\tilde L^R) =  n(\tilde L^R_+)+n(\tilde L^R_-).  
\end{equation}
Thus, it is sufficient to study the kernel and the Morse index of the operators $\tilde L^R_+$ and $\tilde L^R_-$.
\begin{lemma}\label{nLR+}
Let $p>-1$ and the operator $\tilde L^R_+$ be defined by \eqref{decomposition}, then $\ker(\tilde L^R_+)=\{0\}$ and $$n(\tilde L^R_+)=\left\{\begin{array}{c}
     k+1,\,\, \gamma>0 \\
     N-k,\,\, \gamma<0
\end{array}\right.,\,\, k=0, \ldots,\left[\tfrac{N-1}{2}\right].$$
\end{lemma}
\noindent The proof follows from \cite[Proposition 6.1]{AQFF} and \cite[Theorem 3.1]{Kai19}.
Next we evaluate $n(\tilde L^R_-).$ To do that we will use the ideas from \cite{Lop11}.
Denote $\Psi_{k}^{\gamma}=\left(\psi_{k, j}^{\gamma}\right)_{j=1}^{N}, k=0, \ldots,\left[\frac{N-1}{2}\right]$, where 
\begin{equation}\label{hump_1}
\begin{aligned}
\psi_{k, j}^{\gamma}(x)=\left\{\begin{array}{ll}
{\left[p \operatorname{sech}^{2}\Big((p-1) x-a_{k}\Big)\right]^{\frac{1}{2p-2}},} & j=1, \ldots, k ; \\
{\left[p\operatorname{sech}^{2}\Big((p-1)  x+a_{k}\Big)\right]^{\frac{1}{2p-2}},} & j=k+1, \ldots, N,
\end{array}\right. \\
\text { with } a_{k}=\tanh ^{-1}\left(\frac{\gamma}{(N-2k) \sqrt{\omega}}\right), \text { and } \omega>\frac{\gamma^{2}}{(N-2 k)^{2}}.
\end{aligned}
\end{equation} 
We introduce the self-adjoint operator 
\begin{equation}\label{L_e}
\tilde L_\varepsilon=-\frac{d^2}{dx^2}+1-\varepsilon (\Psi_k^\gamma)^{2p-2}, \quad \dom(\tilde L_\varepsilon)=\dom(\Delta_{\frac{\gamma}{\sqrt{\omega}}}),  
\end{equation}
and prove the following technical lemma.
\begin{lemma}\label{Lemma_L_e}
Let the  operator $\tilde L_\varepsilon$ be defined by \eqref{L_e}. Set $$L_{k}^{2}(\Gamma)=\left\{(v_{e})_{e=1}^{N} \in L^{2}(\Gamma): v_{1}(x)=\ldots=v_{k}(x), v_{k+1}(x)=\ldots=v_{N}(x)\right\}.$$ Then the following assertions hold.

$(i)$\, For $\varepsilon<1$ the operator $\tilde L_\varepsilon$  is positive definite.

$(ii)$\, Let $1<\varepsilon< 2p-1$, then the operator $\tilde L_\varepsilon$ is  invertible and 
\begin{itemize}
    \item[$a)$] for $k=0$ in  $L^2(\Gamma)$ we get $n(\tilde L_\varepsilon)=1$ as $\gamma>0$, and $ n(\tilde L_\varepsilon)\geq 1$  as $\gamma<0$;
    \item[$b)$] for $k\geq 1$ in  $L^2_k(\Gamma)$ we get $n(\tilde L_\varepsilon)=1$  as $\gamma<0$, and $n(\tilde L_\varepsilon)\geq 1$  as $\gamma>0$. 
\end{itemize}
\end{lemma}
\begin{proof}
Let $\varepsilon=1$, then it is easily seen that $\ker(\tilde L_1)=\Span\{\Psi_k^\gamma\}$  and $\tilde L_1\geq 0$.
Then assertion $(i)$ is trivial observing that for $\varepsilon<1$ one gets $\tilde L_\varepsilon> \tilde L_1$. 
 
Notice that $\sigma_{\ess}(\tilde L_\varepsilon)=[1, \infty)$, and  for $\varepsilon>1$ we get 
$$(\tilde L_\varepsilon \Psi^\gamma_k, \Psi^\gamma_k)_2=-(\varepsilon-1)\|\Psi^\gamma_k\|_{2p}^{2p}<0,$$ then the first eigenvalue of $\tilde L_\varepsilon$ is negative. In particular, by Min-Max theorem, 
the discrete  spectrum of $\tilde L_\varepsilon$ moves  to the left when $\varepsilon$ increases.
Let $\varepsilon=2p-1$, then by \cite[Proposition 4]{AngGol18a} and \cite[Proposition 3.17]{AngGol18},   we get:

$\bullet$\, $n(\tilde{L}_{2p-1})=1$ for  $k=0, \gamma>0$ in $L^2(\Gamma)$ and for $k\geq 1, \gamma<0$ in $L^2_k(\Gamma)$.

$\bullet$\, $n(\tilde{L}_{2p-1})=2$ for  $k=0, \gamma<0$  and for $k\geq 1, \gamma>0$ in $L^2_k(\Gamma)$.

\noindent Using analyticity  of the family $\tilde L_\varepsilon$, we conclude $a)$ and $b)$  for $1<\varepsilon<2p-1$.
\end{proof}
\begin{lemma}\label{nLR-}
Let the  operator $\tilde L_-^R$ be defined by \eqref{decomposition}. Then the following assertions hold.

$(i)$\, For $b>p-1$ the operator   $\tilde L_-^R$  is positive definite.

$(ii)$\, Let $0<b< p-1$, then the operator $\tilde L_-^R$ is  invertible and 
\begin{itemize}
    \item[$a)$] for $k=0$ in  $L^2(\Gamma)$ we get $n(\tilde L_-^R)=1$ as $\gamma>0$, and $ n(\tilde L_-^R)\geq 1$  as $\gamma<0$;
    \item[$b)$] for $k\geq 1$  in  $L^2_k(\Gamma)$ we get $n(\tilde L_-^R)=1$  as $\gamma<0$, and $n(\tilde L_-^R)\geq 1$  as $\gamma>0$. 
\end{itemize}
\end{lemma}
\begin{proof}
Put $\varepsilon=\frac{2p-1-b}{b+1}$, then $1<\varepsilon<2p-1$ for $0<b<p-1$ and $\varepsilon<1$ for $b>p-1$.
Let $\lambda\in\sigma_p(\tilde L_-^R)$ and  $f(x)\in \dom(\Delta_\gamma)$ be the corresponding eigenvector: $\tilde L_-^R f=\lambda f$. Define $h(x)=f(\frac{x}{\sqrt{\omega}})$, then $h(x)\in \dom(\Delta_{\frac{\gamma}{\sqrt{\omega}}})$ and $\tilde L_\varepsilon h=\frac{\lambda}{\sqrt{\omega}}h.$
\end{proof}
\begin{proof}[Proof of Proposition \ref{count_n}]
The first assertion follows from formula \eqref{decomp_n} and Lemma \ref{nLR+}.  The assertions $a), b)$ follow from Lemma \ref{nLR-}.
\end{proof}
\begin{proof}[Proof of Theorem \ref{main_1}]
 As in the proof of Theorem \ref{spec_inst_1_hump}, we can show that  there exist $\lambda>0$ and $(w_{1}, w_{2}) \in \dom(\Delta_{\gamma})\times \dom(\Delta_{\gamma})$ such that
\begin{equation}\label{lambda}
\left(\begin{array}{cc}
0 & \tilde L^I \\
-\tilde L^R & 0
\end{array}\right)\left(\begin{array}{l}
w_{1} \\
w_{2}
\end{array}\right)=\lambda\left(\begin{array}{c}
w_{1} \\
w_{2}
\end{array}\right).
\end{equation}
As before, we rely on 
\cite[Theorem 1.2]{Grill88} which  states that the number $I(\tilde L^R, \tilde L^I)$ of positive $\lambda$ satisfying \eqref{lambda} is estimated by
\begin{equation}\label{estim_lambda}
n(P\tilde L^R)-n(P(\tilde L^I)^{-1}) \leq I(\tilde L^R, \tilde L^I),
\end{equation}
 where $P$  is the orthogonal projection onto $$\ker(\tilde L^I)^{\perp}=\left\{(v_1,v_2)\in L^{2}(\Gamma)\times L^{2}(\Gamma) :\big((v_1,v_2), (\Phi^\gamma_k,0)\big)_{2}=\big((v_1,v_2), (0,\Phi^\gamma_k)\big)_{2}=0\right\}$$ (here we assume that $L^{2}(\Gamma)$ is endowed with the usual complex inner product).
Using positivity of $\tilde L^I$ and the estimates of $n(\tilde L^R)$ from Proposition \ref{count_n}, formula \eqref{estim_lambda} yields:

$\bullet$\, for $b>-1$: $I(\tilde L^R, \tilde L^I)\geq \left\{\begin{array}{c}
      k-1, \quad \gamma>0 \\
     N-k-2, \quad \gamma<0
\end{array}\right.$

$\bullet$\, for $0<b<p-1$: $I(\tilde L^R, \tilde L^I)\geq \left\{\begin{array}{c}
      k, \quad \gamma>0 \\
     N-k-1, \quad \gamma<0.
\end{array}\right.$

As in the previous subsection, the orbital instability part follows from the fact that the nonlinearity $F_p(u,v)= \left(|u|^{2p-2} u+b|v|^{p}|u|^{p-2} u, |v|^{2p-2} v+b|u|^{p}|v|^{p-2} v\right)$ is of $C^{2}$ class for $p>3$.
\end{proof}
\begin{proof}[Proof of Theorem \ref{main_2}]
Our aim is to use \cite[Stability Theorem and Instablity Theorem]{GrilSha90}. Observe that \cite[\textit{Assumption 2}]{GrilSha90} follows by the Implicit Function Theorem.
Indeed, consider the mapping $F((\omega_1, \omega_2), (u,v)): \mathbb{R}_+^2\times X_{\mathbb{R}} \rightarrow X^*_{\mathbb{R}}$, where $F((\omega_1, \omega_2), (u,v))=S'_{\omega_1,\omega_2}(u,v)$, and $X_{\mathbb{R}}$ is the  Banach space consisting of real-valued functions from $X$. Observe that $\ker(S''_{\omega,\omega}(\Phi^\gamma_0, \Phi^\gamma_0))=\Span\{(0,0)\}$ in $X_{\mathbb{R}}$ since in  $X$ we have $\ker(S''_{\omega,\omega}(\Phi^\gamma_0, \Phi^\gamma_0))=\Span\{(i\Phi^\gamma_0, i\Phi^\gamma_0)\}$.
Observe that $S''_{\omega,\omega}(\Phi^\gamma_0, \Phi^\gamma_0)$  satisfies \cite[Condition (G)]{Stu08}. Indeed, for $\vec{v}=(v_1, v_2)\in X$ 
\begin{equation*}
\begin{split}
 &\langle S''_{\omega,\omega}(\Phi^\gamma_0, \Phi^\gamma_0)\vec{v}, \vec{v}\rangle_{X^*\times X}=\int_{\Gamma}\Big\{|v'_1|^2+|v'_2|^2+\omega\Big(|v_1|^2+|v_2|^2\Big)\Big\}dx
 \\& -\int_{\Gamma}\Big\{(1+b)(\Phi_0^\gamma)^{2p-2}\Big(|v_1|^2+|v_2|^2\Big)+\big(2p-2+b(p-2)\big)(\Phi_0^\gamma)^{2p-2}\Big((\re v_1)^2+(\re v_2)^2\Big)\Big\}dx\\
 & -2bp\int_{\Gamma} (\Phi_0^\gamma)^{2p-2}\re v_1 \re v_2\Big\}dx-\gamma\Big(|v_{11}(0)|^2+|v_{21}(0)|^2\Big)\\&
 \geq \min\{\frac{1}{2}, \omega\}\|\vec{v}\|^2_X-\Big((1+b)(2p-1)M^{2p-2}+\frac{2\gamma^2}{N^2}\Big)\|\vec{v}\|^2_2,\quad M=\|\Phi_0^\gamma\|_\infty.
 \end{split}    
\end{equation*}

Thus, by \cite[Lemma 5.4]{Stu08}, we conclude $\sigma_{\ess}(\mathcal{R}^{-1}S''_{\omega,\omega}(\Phi^\gamma_0, \Phi^\gamma_0))=\sigma_{\ess}(\tilde L)=[\omega, \infty)$ since $\tilde L=\tilde L^R+i \tilde L^I$ and $\sigma_{\ess}(\tilde L^R)=\sigma_{\ess}(\tilde L^I)=[\omega, \infty)$, where  $\mathcal{R}: X \rightarrow X^{*}$ is the Riesz isomorphism. Hence $S''_{\omega,\omega}(\Phi^\gamma_0, \Phi^\gamma_0):X_\mathbb{R}\to X_\mathbb{R}^*$ is invertible and bounded.  Then, by the Implicit Function Theorem,   there is an open neighborhood $\Omega$ of $(\omega, \omega)$ and a unique function
$(\Phi_1(\omega_1, \omega_2), \Phi_2(\omega_1, \omega_2)): \Omega \rightarrow X_{\mathbb{R}}$ such that $S'_{\omega_1,\omega_2}(\Phi_1(\omega_1, \omega_2), \Phi_2(\omega_1, \omega_2))=S'_{\omega,\omega}(\Phi_0^\gamma, \Phi_0^\gamma)=0$. 

Next, let us check  that  \cite[\textit{Assumption} 3]{GrilSha90} holds.
By \eqref{decomp_ker}, \eqref{decomp_n},  Lemmas \ref{nLR+} and \ref{nLR-}, we conclude that $\ker(\tilde L^R)=\{0\}$ for $b\neq p-1$,  and $n(\tilde L^R)< \infty$,  then, by \eqref{4.13}  and \cite[Lemma 5.4]{Stu08}, we obtain $\ker(\mathcal{R}^{-1}S''_{\omega,\omega}(\Phi^\gamma_0, \Phi^\gamma_0))=\ker(\tilde L)=\Span\{i(\Phi^\gamma_0,0), i(0, \Phi^\gamma_0)\}$
 and $n(\mathcal{R}^{-1}S''_{\omega,\omega}(\Phi^\gamma_0, \Phi^\gamma_0))=n(\tilde L^R)<\infty$.   
 
 Let $d''(\omega)$ be  the Hessian of $S_{\omega_1, \omega_2}(\Phi_1(\omega_1, \omega_2), \Phi_2(\omega_1, \omega_2))$ at $(\omega, \omega)$.
 It is easily seen that $$d''(\omega)=\frac{1}{2}\left(\begin{array}{cc}
   \partial_{\omega_1}\|\Phi_1(\omega_1,\omega_2)\|^2_2  & \partial_{\omega_2}\|\Phi_1(\omega_1,\omega_2)\|^2_2  \\
   \partial_{\omega_1}\|\Phi_2(\omega_1,\omega_2)\|^2_2  & \partial_{\omega_2}\|\Phi_2(\omega_1,\omega_2)\|^2_2
\end{array}\right)\Big|_{(\omega_1, \omega_2)=(\omega, \omega)}. $$
Denote by $p(d''(\omega))$ the number of positive eigenvalues of $d''(\omega)$.
By \cite[Stability Theorem and Instability Theorem]{GrilSha90}, we conclude:
 
 $\bullet$\, if  $n(\tilde L)=p(d''(\omega))$, then $e^{i\omega t}(\Phi^\gamma_0, \Phi^\gamma_0)$ is orbitally stable;
 
 $\bullet$\, if $n(\tilde L)-p(d''(\omega))$ is odd, then $e^{i\omega t}(\Phi^\gamma_0, \Phi^\gamma_0)$ is spectrally unstable.

\textit{Step 1.} Note that $\partial_{\omega_j}\|\Phi_i(\omega_1,\omega_2)\|^2_2=\int_\Gamma \Phi_i(\omega_1,\omega_2) \partial_{\omega_j}\Phi_i(\omega_1,\omega_2)dx$,  $i,j=1,2$.
Differentiating $S'_{\omega, \omega}(\Phi_1(\omega_1, \omega_2), \Phi_2(\omega_1, \omega_2))$ with respect to $\omega_1$, using the chain rule, and denoting $(h_1, h_2):=(\partial_{\omega_1}\Phi_1, \partial_{\omega_1}\Phi_2)|_{(\omega_1, \omega_2)=(\omega, \omega)}$, we obtain for $v\in H^1(\Gamma)$
\begin{align*}
       &\int_{\Gamma}\Big\{h'_1\overline{v}'+\omega h_1 \overline{v}-(2p-1+b)(\Phi^\gamma_0)^{2p-2}h_1\overline{v}-bp(\Phi^\gamma_0)^{p-2}h_2\overline{v}\Big\}dx-\gamma h_{11}(0)\overline{v}_1(0)=-\int_{\Gamma} \Phi^\gamma_0\overline{v}dx    \\
      &\int_{\Gamma}\Big\{h'_2\overline{v}'+\omega h_2 \overline{v}-(2p-1+b(p-1))(\Phi^\gamma_0)^{2p-2}h_2\overline{v}-bp(\Phi^\gamma_0)^{p-2}h_1\overline{v}\Big\}dx-\gamma h_{21}(0)\overline{v}_1(0)=0.  
  \end{align*}
  Analogous system can be obtained when differentiating with $\omega_2$.
  Then setting $h_+=h_1+h_2$ and $h_-=h_1-h_2, $ we have  
\begin{equation}\label{syst_quad_form}
  \begin{split}
    \left\{\begin{array}{c}
      \mathfrak{t}_+(h_+, v)=-(\Phi^\gamma_0, v)_2,    \\
     \mathfrak{t}_-(h_-, v)=-(\Phi^\gamma_0, v)_2,  
    \end{array}\right.
    \end{split}
    \end{equation}
    where $\mathfrak{t}_+$ and $\mathfrak{t}_-$ are bilinear  forms associated with the operators $\tilde L^R_+$ and $\tilde L^R_-$ defined by \eqref{decomposition}. Hence 
    \begin{equation}\label{H}d''(\omega)=\frac{1}{4}\left(\begin{array}{cc}
   (\Phi^\gamma_0, h_+ + h_-)_2 &  (\Phi^\gamma_0, h_+ - h_-)_2   \\
 (\Phi^\gamma_0, h_+ - h_-)_2   & (\Phi^\gamma_0, h_+ + h_-)_2  \end{array}\right). \end{equation}
 Observe that $\det(d''(\omega))=\frac{(\Phi^\gamma_0, h_+)_2(\Phi^\gamma_0, h_-)_2}{2}$ and $\mathrm{trace}(d''(\omega))=\frac{(\Phi^\gamma_0, h_++h_1)_2}{2}$.

 \textit{Step 2.} Below we will prove:
 
 $a)$\, if $b>-1$, then  $(\Phi^\gamma_0, h_+)_2>0$ for $2\leq p<3$, and $(\Phi^\gamma_0, h_+)_2<0$ for $p>3$;
 
 $b)$\, if     $b>p-1$, then  $(\Phi^\gamma_0, h_-)_2<0$, and    if  $0<b<p-1$, then $(\Phi^\gamma_0, h_-)_2>0$.
 
 Firstly, we prove $a)$. Since $\Phi^\gamma_0$
 satisfies \eqref{profile_multi} and $v\in H^1(\Gamma)$ in \eqref{syst_quad_form} is arbitrary, we conclude that $h_+=\partial_\omega \Phi^\gamma_0$. Then $(\Phi^\gamma_0, h_+)_2=(\Phi^\gamma_0, \partial_\omega \Phi^\gamma_0)_2=\frac{1}{2}\partial_\omega(\Phi^\gamma_0, \Phi^\gamma_0)_2$.
 Moreover, using \eqref{hump_1}, we have 
$$
\begin{aligned}
\int_{\Gamma} (\Phi^\gamma_0)^2(x) d x &=\left(\frac{\omega}{b+1}\right)^{1 /(p-1)} \int_{\Gamma} (\Psi_{0}^\gamma)^{2}(\sqrt{\omega} x) d x \\
&=\frac{\omega^{(3-p) /(2 p-2)}}{(b+1)^{1 /(p-1)}} \int_{\Gamma} (\Psi_{0}^\gamma)^{2}(x) dx,
\end{aligned}
$$
and then
$$
\partial_\omega(\Phi^\gamma_0, \Phi^\gamma_0)_2=\frac{(3-p)}{(2 p-2)} \frac{\omega^{(5-3 p) /(2 p-2)}}{(b+1)^{1 /(p-1)}} \int_{\Gamma} (\Psi^\gamma_{0})^{2}(x) d x
$$
which is positive for $2 \leq p<3$ and negative for $p>3$. This proves $a)$.

Secondly, we  prove $b)$ for $b>p-1$.
Notice that the second line of \eqref{syst_quad_form} is equivalent to 
\begin{equation}\label{1}
\begin{split}
&\int_{\Gamma}\Big\{h'_-(x)\overline{v}'(x)+\omega h_-(x)\overline{v}(x)-\frac{\omega(2p-b-1)}{b+1}\big(\Psi^\gamma_0(\sqrt{\omega}x)\big)^{2p-2}h_-(x)\overline{v}(x)\Big\}dx\\&-\gamma h_{-,1}(0)\overline{v}_1(0)=-\left(\frac{\omega}{b+1}\right)^{1/2(p-1)}\int_{\Gamma}\Psi_0^\gamma(\sqrt{\omega}x)\overline{v}(x)dx.  \end{split} 
\end{equation}
Denoting $f(\sqrt{\omega}x)=v(x), s(\sqrt{\omega}x)=h_-(x)$, from \eqref{1} we get 
\begin{equation}\label{2}
\begin{split}
&\int_{\Gamma}\Big\{s'(y)\overline{f}'(y)+ s(y)\overline{f}(y)-\frac{(2p-b-1)}{b+1}\big(\Psi^\gamma_0(y)\big)^{2p-2}s(y)\overline{f}(y)\Big\}dy-\frac{\gamma}{\sqrt{\omega}} s(0)\overline{f}_1(0)\\&=-\frac{\omega^{(3-2p)/(2p-2)}}{(b+1)^{1/(2p-2)}}\int_{\Gamma}\Psi_0^\gamma(y)\overline{f}(y)dy.  \end{split} 
\end{equation}
Moreover, 
$$\int_{\Gamma}h_-(x)\Psi^\gamma_0(\sqrt{\omega}x)dx=\frac{1}{\sqrt{\omega}}\int_{\Gamma}h_-(\frac{y}{\sqrt{\omega}})\Psi^\gamma_0(y)dy=\frac{1}{\sqrt{\omega}}\int_{\Gamma}s(y)\Psi^\gamma_0(y)dy,$$ then $(h_-, \Phi^\gamma_0)_2$ and  $(s, \Psi^\gamma_0)_2$ have the same sign.
Observe that \eqref{2} can be rewritten as $\mathfrak{t}_\varepsilon(s,f)=-C(\omega)(\Psi^\gamma_0,f)_2$, where $\mathfrak{t}_\varepsilon$ is the bilinear  form associated with self-adjoint operator $\tilde L_\varepsilon$ given by \eqref{L_e} (for $\varepsilon=\frac{2p-b-1}{b+1}$) and $C(\omega)=\frac{\omega^{(3-2p)/(2p-2)}}{(b+1)^{1/(2p-2)}}$.  Set $f=s$, then by Lemma \ref{Lemma_L_e}-$(i)$, we conclude that $0<\mathfrak{t}_\varepsilon(s,s)=-C(\omega)(\Psi^\gamma_0, s)_2$  for $b>p-1$. 

Thirdly, we prove $b)$ for $0<b<p-1$. Notice that from \eqref{2},  by the Representation Theorem  \cite[Chapter VI, Theorem 2.1]{Kat66}, we conclude 
\begin{equation}\label{3} s\in \dom(\tilde L_\varepsilon), \qquad \tilde L_\varepsilon s=-C(\omega)\Psi^\gamma_0.\end{equation}
Since $\tilde L_\varepsilon$ is invertible for $1<\varepsilon=\frac{2p-b-1}{b+1}<2p-1$ (or equivalently for  $0<b<p-1$) and holomorphic in $\varepsilon$, then the solution $s=s(\varepsilon)$ of \eqref{3} is smooth in $\varepsilon$. Differentiating $\mathfrak{t}_{\varepsilon}(s(\varepsilon), \cdot)$ with $\varepsilon$, from \eqref{3} we obtain 
$$\mathfrak{t}_{\varepsilon}\big(\partial_\varepsilon s(\varepsilon),f\big)=\big(s(\varepsilon)(\Psi_0^\gamma)^{2p-2}, f\big)_2.$$
Let $f=s(\varepsilon),$ then again, by \eqref{3},
\begin{align*}
 &0<\int_{\Gamma} |s(\varepsilon)|^2(\Psi^\gamma_0(x))^{2p-2}dx=\mathfrak{t}_\varepsilon\big(\partial_\varepsilon s(\varepsilon), s(\varepsilon)\big)=\big(\partial_\varepsilon s(\varepsilon), \tilde L_\varepsilon s(\varepsilon)\big)_2\\&=-C(\omega)\big(\partial_\varepsilon s(\varepsilon), \Psi^\gamma_0\big)_2= -C(\omega)\partial_\varepsilon\big(s(\varepsilon), \Psi^\gamma_0\big)_2. \end{align*}
Hence $(s(\varepsilon), \Psi^\gamma_0)_2$ is decreasing.

It is easily seen that  $v(\omega,x)=\omega^{1/(2p-2)}\Psi^\gamma_0(\sqrt{\omega}x)$ satisfies 
$$-\Delta_\gamma v(\omega,x)+\omega v(\omega,x)-v(\omega,x)^{2p-1}=0.$$
Further, from the above one concludes that $u(x)=\partial_\omega v(\omega,x)|_{\omega=1}$ is the solution to
$$-\Delta_\gamma u(x)+u(x)-(2p-1)\big(\Psi_0^\gamma(x)\big)^{2p-2}u(x)=-\Psi_0^\gamma(x).$$
Observe that $\partial_\omega v(\omega,x)\in\dom(\Delta_\gamma)$. Indeed,  define $$\hat{S}_\omega(\psi)=\frac{1}{2}\left\{\|\psi'\|_2^2-\gamma|\psi_{1}(0)|^2+\omega\|\psi\|_2^2\right\}-\frac{1}{2p}\|\psi\|_{2p}^{2p},$$ then $\hat{S}'_\omega(v(\omega, x))=0$ and therefore $\left\langle\hat{S}''_\omega(v(\omega, x))\partial_\omega v(\omega, x), g\right\rangle_{(H^1)^*\times H^1}=0, \, g\in H^1(\Gamma)$, which yields 
\begin{align*}
&\int_{\Gamma}(\partial_\omega v(\omega, x))'\overline{g}(x)dx-\gamma\partial_\omega v_1(\omega, 0)\overline{g}_1(0)\\&=\int_{\Gamma}\left\{-\omega\partial_\omega v(\omega, x)\overline{g}(x)+\Big(\partial_\omega v(\omega, x)+(2p-2)\re(\partial_\omega v(\omega, x))\Big)\omega(\Psi_0^\gamma(\sqrt{\omega}x))^{2p-2} \overline{g}(x)\right\}dx.    
\end{align*}
By the Representation Theorem, $\partial_\omega v(\omega,x)\in\dom(\Delta_\gamma)$.

It is easily seen that $\partial_\omega v(\omega,x)|_{\omega=1}C(\omega)=s(2p-1),$ where $s(2p-1)$ is the solution to \eqref{3}  for $\varepsilon=2p-1$. By the continuity of the scalar product, we conclude
\begin{align*}
&\big(s(2p-1), \Psi^\gamma_0\big)_2=\big(\partial_\omega v(\omega,x)|_{\omega=1}C(\omega), \Psi^\gamma_0\big)_2=C(\omega)\frac{1}{2}\partial_\omega \|v(\omega,x)\|^2|_{\omega=1}\\&=\frac{1}{2}C(\omega)\partial_\omega\Big(\omega^{(3-p)/(2p-2)}\int_{\Gamma}(\Psi^\gamma_0)^2dx\Big)\Big|_{\omega=1}>0\quad\text{for}\quad p<3.\end{align*}
Finally, since $(s(\varepsilon), \Psi^\gamma_0)_2$ is decreasing on the interval $(1, 2p-1)$, we get $(s(\varepsilon), \Psi^\gamma_0)_2>0$  for $\varepsilon\in (1, 2p-1)$.

\textit{Step 3.}  Let $\gamma>0, k=0$.  By \eqref{decomp_n},  Lemmas \ref{nLR+} and \ref{nLR-}, we have    $n(\tilde L^R)=1$ for $b>p-1$, and    $n(\tilde L^R)=2$ for $0<b<p-1$. 

Let $2\leq p< 3$, then, by \textit{Step 2} and \eqref{H}, we conclude that  $p(d''(\omega))=1$  for $b>p-1$, and   $p(d''(\omega))=2$ for $0<b<p-1$.
Let $p>3$ and $b>p-1$, then $p(d''(\omega))=0$. 

Finally, for $2\leq p<3$ and $b\neq p-1$, we get $n(\tilde L)=p(d''(\omega))$, therefore,  $e^{i\omega t}(\Phi^\gamma_0,\Phi^\gamma_0)$ is orbitally stable.  

For  $p>3$ and $b>p-1$, we obtain $n(\tilde L)- p(d''(\omega))=1$, then, by\cite[Theorem 5.1]{GrilSha90}, there exists  positive $\lambda$ satisfying  \eqref{lambda}, and orbital instability follows by $C^2$ regularity of the data-solution mapping associated with \eqref{syst_particular}.

\end{proof}
\subsection{Stability of the standing wave generated by 2D rotation}\label{subsec_3.3}
Observe that for $p=2$ and $b=1$, in addition to the symmetry $T_0(\theta_0)(u,v)=e^{i\theta_0}(u,v), \theta_0\in\mathbb{R}$,   system  \eqref{syst_particular} is also invariant under 2D rotation $T_1(\theta_1)$:  $$
\binom{u(t, x)}{v(t, x)} \overset{T_1(\theta_1)}{\mapsto}\left(\begin{array}{rr}
\cos \theta_1 & \sin \theta_1 \\
-\sin \theta_1 & \cos \theta_1
\end{array}\right)\binom{u(t, x)}{v(t, x)}, \quad \theta_1 \in \mathbb{R}.
$$ 
Below we study orbital stability of the standing wave related to both of these symmetries: 
$$(u(t, x),v(t, x))=T_0(\omega_0 t)T_1(\omega_1 t)\binom{\Phi_1(x)}{\Phi_2(x)}, \quad \omega_0,\omega_1\in \mathbb{R}^+.$$ It is easily seen that $T'_1(0)=\left(\begin{array}{rr}
0 & 1 \\
-1 & 0 
\end{array}\right). $
Assume that the profile of the standing wave has the form $(\Phi_1, \Phi_2)=\Phi_{\omega_0,\omega_1}(x)\cdot(v_1,v_2)$, where  $(v_1,v_2)=(\frac{i}{\sqrt{2}}, \frac{1}{\sqrt{2}})$ is the unit eigenvector of $iT'_1(0)$ corresponding to the eigenvalue $1$. Then the profile $\Phi_{\omega_0,\omega_1}(x)$ satisfies  the stationary equation
\begin{equation}\label{sol_2D}
    -\Delta_\gamma \Phi+(\omega_0-\omega_1)\Phi-|\Phi|^2\Phi=0.
\end{equation}
The description of solutions is given by \eqref{hump} (with $\omega$ substituted by $\omega_0-\omega_1$).
Below we will prove the following theorem.
\begin{theorem}\label{main_2D}
Let $\omega_0-\omega_1>\frac{\gamma^2}{N^2}$  and $\Phi^\gamma_k$ be defined by \eqref{hump}, then for $k=0, \gamma>0$ the standing wave $T_0(\omega_0 t)T_1(\omega_1 t)\binom{i\Phi_k^\gamma/\sqrt{2}}{\Phi^\gamma_k/\sqrt{2}}$ is orbitally stable.
Moreover, if either $\gamma>0$ and $k$ is odd, or $\gamma<0$ and $N-k-1$ is odd, then $T_0(\omega_0 t)T_1(\omega_1 t)\binom{i\Phi_k^\gamma/\sqrt{2}}{\Phi^\gamma_k/\sqrt{2}}$ is spectrally unstable.
\end{theorem}
\noindent In the above theorem the stability is understood in the sense of the $U(1)$-symmetry. This is related to the fact that the proof uses \cite[Stability Theorem and Instablity Theorem]{GrilSha90} where the stability/instability  results are stated for a  centralizer subgroup. %It is easily seen that centralizer of $\{T_0(\theta_0)T_1(\theta_1)\}$  coincides with $T_0(\cdot)$. 

Conserved functional generated by $T'_1(0)$ is given by 
$$Q_1(u,v)=\im \int_{\Gamma}u\overline{v}dx,$$ and the standing wave $T_0(\omega_0 t)T_1(\omega_1 t)\binom{i\Phi^\gamma_k/\sqrt{2}}{\Phi^\gamma_k/\sqrt{2}}$ is a critical point of the  functional 
\begin{equation*}
    S_{\omega_0,\omega_1}(u,v)=\frac{1}{2}\left\{E(u,v)+\omega_0(\|u\|_2^2+\|v\|_2^2)\right\}+\omega_1Q_1(u,v).
\end{equation*}
As in the previous subsections,  we associate with the  bilinear form 

$\left\langle S_{\omega_0,\omega_1}^{\prime \prime}\left(\frac{i\Phi^\gamma_k}{\sqrt{2}}, \frac{\Phi^\gamma_k}{\sqrt{2}}\right)\left(h_{1}, h_{2}\right),\left(z_{1}, z_{2}\right)\right\rangle_{X^{*} \times X}$ the self-adjoint in $L^2(\Gamma)\times L^2(\Gamma)$ operator  $\tilde{L}$. Let $h_j=h_j^R+ih_j^I, j=1,2$.  Substituting complex-valued vector function $\vec{h}= (h_1, h_2)$ by the corresponding quadruplet of
real-valued functions $(h_1^R, h_1^I, h_2^R, h_2^I)$ and $\tilde L\vec{h}=(f,g)=(f^R+if^I, g^R+ig^I)$  by the quadruplet $(f^R, f^I, g^R, g^I)$, we can express  the action of $\tilde L$ as 
\begin{equation}\label{L_2D}\left(\begin{array}{c}
     f^R  \\
     f^I \\
     g^R \\
     g^I
     \end{array}\right)=\left(\begin{array}{cccc}
     L_{\omega_0}-(\Phi^\gamma_k)^2    &  0 & 0 & \omega_1\\
      0 &   L_{\omega_0}-2(\Phi^\gamma_k)^2 & -\omega_1-(\Phi^\gamma_k)^2  & 0\\
       0 &  -\omega_1-(\Phi^\gamma_k)^2 &  L_{\omega_0}-2(\Phi^\gamma_k)^2 &    0\\
       \omega_1 & 0&0& L_{\omega_0}-(\Phi^\gamma_k)^2      \end{array}\right)\left(\begin{array}{c}
     h^R_1  \\
     h^I_1 \\
     h^R_2 \\
     h^I_2
     \end{array}\right),\end{equation}
     where $L_{\omega_0}=-\Delta_\gamma+\omega_0$ with $\dom(L_{\omega_0})=\dom(\Delta_\gamma)$.
     Below we characterize the spectral properties of the operator $\tilde L$.
     \begin{lemma}\label{spec_L_2D}
     Let $\omega_0-\omega_1>\frac{\gamma^2}{N^2}$ and $\Phi_{\omega_0,\omega_1}$ be the solution to \eqref{sol_2D} given by \eqref{hump} (with $\omega$ substituted by $\omega_0-\omega_1$). Let also operator $\tilde L$ be defined by \eqref{L_2D}.Then the following assertions hold:
     
     $(i)$ $\ker(\tilde L)=\Span\{(\Phi^\gamma_k, -i \Phi^\gamma_k)\} $;
     
     $(ii)$ $n(\tilde L)=\left\{\begin{array}{c}
          k+1, \, \gamma>0 \\
          N-k, \, \gamma<0
     \end{array}\right.; $
     
     $(iii)$ $\sigma_{\ess}(\tilde L)>0$.
     \end{lemma}
     \begin{proof}
 Observe that $(\tilde L-\lambda)\vec{h}=(f,g)$  is equivalent to 
   \begin{equation}\label{system_2D}\left\{\begin{array}{c}        (L_{\omega_0}+\omega_1-(\Phi^\gamma_k)^2-\lambda)(h_1^R+h_2^I)=f^R+g^I  \\
        (L_{\omega_0}-\omega_1-(\Phi^\gamma_k)^2-\lambda)(h_1^R-h_2^I)=f^R-g^I\\
        (L_{\omega_0}-\omega_1-3(\Phi^\gamma_k)^2-\lambda)(h_1^I+h_2^R)=f^I+g^R\\
        (L_{\omega_0}+\omega_1-(\Phi^\gamma_k)^2-\lambda)(h_1^I-h_2^R)=f^I-g^R.
   \end{array}\right.\end{equation}
   By \cite[Proposition 6.1]{AQFF}, the operator $L_{\omega_0}-\omega_1-(\Phi^\gamma_k)^2 $ is nonnegative and $\ker(L_{\omega_0}-\omega_1-(\Phi^\gamma_k)^2)=\Span\{\Phi^\gamma_k\}$. Moreover, $\ker(L_{\omega_0}-\omega_1-3(\Phi^\gamma_k)^2)=\{0\}$ and, by \cite[Theorem 3.1]{Kai19},  $n(L_{\omega_0}-\omega_1-3(\Phi^\gamma_k)^2)=\left\{\begin{array}{c}
          k+1, \, \gamma>0 \\
          N-k, \, \gamma<0
     \end{array}\right..$
   Observing that operator $L_{\omega_0}+\omega_1-(\Phi^\gamma_k)^2$ is positive definite, we arrive at $(i)$, $(ii)$. 
   
   Finally, noticing that   for $\lambda\in\mathbb{R}^+\setminus[\omega_0-\omega_1, \infty)$ all the operators on the left side of \eqref{system_2D} are invertible, we conclude that  $\lambda\in\rho(\tilde L)$ (since $\tilde L-\lambda$ is bijective), and therefore $(iii)$ follows. 
        \end{proof}
        
\begin{proof}[Proof of Theorem \ref{main_2D}]
Without abuse of notation, we will write $S''_{\omega_0,\omega_1}$  instead of $S''_{\omega_0,\omega_1}\left(\frac{i\Phi^\gamma_0}{\sqrt{2}}, \frac{\Phi^\gamma_0}{\sqrt{2}}\right)$.  Observe that $S''_{\omega_0,\omega_1}$  satisfies \cite[Condition (G)]{Stu08}. Indeed, for $\vec{v}=(v_1, v_2)\in X$ 
\begin{equation*}
\begin{split}
 \left\langle S''_{\omega_0,\omega_1}\vec{v}, \vec{v}\right\rangle_{X^*\times X}&=\int_{\Gamma}\Big\{|v'_1|^2+|v'_2|^2+\omega_0\Big(|v_1|^2+|v_2|^2\Big)\Big\}dx
 \\& -\int_{\Gamma}\Big\{(\Phi_0^\gamma)^{2}\Big(|v_1|^2+|v_2|^2+(\re v_1)^2+(\re v_2)^2+2\re v_2\im v_1\Big)\Big\}dx\\
 & -2\omega_1\im\int_{\Gamma} v_1\overline{v_2}dx-\gamma\Big(|v_{11}(0)|^2+|v_{21}(0)|^2\Big)\\&
 \geq \min\{\frac{1}{2}, \omega_0\}\|\vec{v}\|^2_X-\Big(3M^{2}+\frac{2\gamma^2}{N^2}+\omega_1\Big)\|\vec{v}\|^2_2,
 \end{split}    
\end{equation*}
where $M=\|\Phi_0^\gamma\|_\infty$. 
Thus, by \cite[Lemma 5.4]{Stu08} and  Lemma \ref{spec_L_2D},  we conclude
\begin{equation*}\begin{split}\label{spec2D}\sigma_{\ess}(\mathcal{R}^{-1}S''_{\omega_0,\omega_1})=\sigma_{\ess}(\tilde L)&>0,\quad \ker(\mathcal{R}^{-1}S''_{\omega_0,\omega_1})=\Span\{(\Phi^\gamma_k, -i\Phi^\gamma_k)\},\\ & n(\mathcal{R}^{-1}S''_{\omega_0,\omega_1})=n(\tilde L).\end{split} \end{equation*}

Since the centralizer subgroup of the group $\{T_0(\theta_0)T_1(\theta_1): \theta_0, \theta_1\in \mathbb{R}\}$ coincides with $G_{\theta_0, \theta_1}=\{T_0(\theta): \theta\in \mathbb{R}\}$, then according to \cite[Stability Theorem and Instability Theorem]{GrilSha90}, we can analyze stability only in context of the $U(1)$-symmetry.  
It is easily seen that \textit{Assumptions 1-3} of \cite{GrilSha90} are satisfied (in particular, \textit{Assumption 3} is satisfied due to Lemma \ref{spec_L_2D}).   
 
 Let $d''(\omega_0-\omega_1)$ be  the reduced  Hessian of $S_{\omega_0, \omega_1}\left(\frac{i\Phi^\gamma_k}{\sqrt{2}}, \frac{\Phi^\gamma_k}{\sqrt{2}}\right)$. Thus, by \cite[Stability Theorem and Instability Theorem]{GrilSha90}, we conclude:
 
 $\bullet$\, if  $n(\tilde L)=p(d''(\omega_0-\omega_1))$, then $T_0(\omega_0 t)T_1(\omega_1 t)\binom{i\Phi_k^\gamma/\sqrt{2}}{\Phi^\gamma_k/\sqrt{2}}$ is orbitally stable;
 
 $\bullet$\, if $n(\tilde L)-p(d''(\omega_0-\omega_1))$ is odd, then $T_0(\omega_0 t)T_1(\omega_1 t)\binom{i\Phi_k^\gamma/\sqrt{2}}{\Phi^\gamma_k/\sqrt{2}}$ is spectrally unstable.
Using \cite[Proposition 5]{AngGol18a}, we conclude  $d''(\omega_0-\omega_1)=\frac{1}{2}\partial_{\omega_0}\|\Phi^\gamma_k\|^2_2>0.$ Finally, applying Lemma \ref{spec_L_2D}-$(ii)$, we arrive at the result. 
\end{proof}        
\section*{Appendix}

In this Appendix we recall some basic properties of symmetric rearrangements of a measurable function $u:\Gamma\rightarrow\mathbb{C}^{N}$. 
Set
\begin{equation*}
\mu_{u}(s)=\left|\{x\in \Gamma:\left|u(x)\right|\geq s\}\right|\,\,\,\,\,\mbox{and}\,\,\,\,\,\lambda_{u}(t)=\sup\left\{s\,:\mu_{u}(s)>Nt\,\right\}.
\end{equation*}
Here $\{x\in \Gamma:\left|u(x)\right|\geq s\}=\bigcup\limits_{e=1}^N\{x\in I_e: |u_e(x)|\geq s\},$ and all sets in the union are disjoint. 
The symmetric rearrangement $u^{*}$ of $u$ is defined by $u^{*}=\left(u^{*}_{e}\right)^{N}_{e=1}$ with $u^{*}_{e}=\lambda_{u}$ for all $e=1,\ldots,N$. Basic properties of symmetric rearrangements are listed in the next proposition.

\begin{proposition} \label{Pr0p0}
Let $u,v\in H^{1}(\Gamma)$. Then the following assertions hold.\\
$(i)$\,The symmetric rearrangement $u^{*}$ is positive and nonincreasing. Moreover, $\left|u\right|$ and $\left|u^{*}\right|$ are equimeasurable, that is, for every $s>0$
\begin{equation}\label{equim}
\left|\{\left|u\right|\geq s\}\right|=\left|\{u^{*}\geq s\}\right|.
\end{equation}
$(ii)$ $u^{*}\in H^{1}(\Gamma)$, $\left\|u\right\|_{p}=\left\|u^{*}\right\|_{p}$ for all $1\leq p\leq\infty$, and
\begin{equation}\label{deriv_rearrag}
\left\|\left(u^{*}\right)'\right\|_{2}\leq\frac{N}{2}\left\|u'\right\|_{2}.
\end{equation}
$(iii)$ If $u,v$ are nonnegative, then
\begin{equation*}
 (u,v)_{2}\leq(u^{*},v^{*})_{2}.
\end{equation*}
\end{proposition}
\begin{proof} The proof of statements $(i)$ and $(ii)$ is contained in \cite[Proposition A.1, Theorem 6]{AQFF}. It follows from the  Layer cake representation theorem \cite[Theorem 1.13]{ELL} and \eqref{equim} that
\begin{equation*}
\begin{split}
&\int_{\Gamma}u(x)v(x)dx=\sum^{N}_{e=1}\int_{\mathbb{R}^{+}}u_{e}(x)v_{e}(x)dx=\sum^{N}_{e=1}\int_{\mathbb{R}^{+}}\left(\int\limits_0^\infty\chi_{\{u_e\geq t\}}(x)dt\int\limits_0^\infty\chi_{\{v_e\geq s\}}(x)ds\right)dx\\
&=\sum^{N}_{e=1}\int\limits_0^\infty\int\limits_0^\infty\left|\left\{u_e\geq s\right\}\cap\left\{v_e\geq t\right\}\right|dsdt=\int^{\infty}_{0}\int^{\infty}_{0}\left|\left\{u\geq s\right\}\cap\left\{v\geq t\right\}\right|\,dsdt\\
&\leq\int^{\infty}_{0}\int^{\infty}_{0}\min\left\{\left|\left\{u\geq s\right\}\right|\,;\,\left|\left\{v\geq t\right\}\right|\right\}\,dsdt
=\int^{\infty}_{0}\int^{\infty}_{0}\min\left\{\left|\left\{u^{*}\geq s\right\}\right|\,;\,\left|\left\{v^{*}\geq t\right\}\right|\right\}\,dsdt\\
&=\int^{\infty}_{0}\int^{\infty}_{0}\left|\left\{u^{*}\geq s\right\}\cap\left\{v^{*}\geq t\right\}\right|\,dsdt=\int_{\Gamma}u^{*}(x)v^{*}(x)dx.
\end{split}
\end{equation*}
Hence $(iii)$ holds.
\end{proof}

%\bibliographystyle{plain}
%\bibliography{bibliografia} 

\end{document}